\documentclass[12pt]{amsart}
\usepackage{hyperref}

\hypersetup{
  colorlinks= true, 
  urlcolor   = red, 
  linkcolor  = blue, 
  citecolor  = blue 
}

\newcommand{\modd}[1]{\,(\text{mod}\, {#1})}
\def\Re{\text{Re}\,}
\def\Im{\text{Im}\,}

\newtheorem{theorem}{Theorem}
\newtheorem{lemma}{Lemma}
\newtheorem*{corollary}{Corollary}
\theoremstyle{definition}
\newtheorem{example}{Example}

\makeatletter
\def\specialsection{\@startsection{section}{1}%
  \z@{\linespacing\@plus\linespacing}{.5\linespacing}%
  {\normalfont}}
\def\section{\@startsection{section}{1}%
  \z@{.7\linespacing\@plus\linespacing}{.5\linespacing}%
  {\normalfont\scshape}}
\makeatother

\title{Harmonic Analysis on the Positive Rationals I:  Basic Results}

\author{Peter D. T. A. Elliott}
\address{Department of Mathematics, University of Colorado Boulder, Boulder, Colorado 80309-0395 USA}
\email{pdtae@euclid.colorado.edu}

\author{Jonathan Kish}
\address{Department of Mathematics, University of Colorado Boulder, Boulder, Colorado 80309-0526 USA}
\email{jonathan.kish@colorado.edu}

\begin{document}
\parskip10pt
\parindent15pt
\baselineskip15pt

\maketitle

\section{Introduction}

A complex-valued function, $g$, is \emph{arithmetic} if it is defined on the positive integers.  It is \emph{multiplicative} if it satisfies $g(ab) = g(a)g(b)$ on mutually prime integers $a$, $b$; \emph{completely multiplicative} if it further satisfies $g(p^k) = g(p)^k$ on prime-powers; \emph{exponentially multiplicative} if $g(p^k) = g(p)^k/k!$.

This study is in two parts.  The present paper contains a complete proof of


\makeatletter{}

\begin{theorem}\label{JP_elliott_03_thm_01}
Let $\alpha$, $\gamma$, $x$ be real numbers and $D$ an integer satisfying $0 < \alpha < 1$, $0 < \gamma < 1$, $1 \le D \le x$.  Let $g$ be a multiplicative function with values in the complex unit disc.

Then there are nonprincipal Dirichlet characters $\chi_j \modd{D}$, their number bounded in terms of $\alpha$ alone, such that
\begin{align}
\sum\limits_{\substack{n \le y \\ n \equiv a \modd{D}}} g(n) = \frac{1}{\varphi(D)}& \sum\limits_{\substack{n \le y \\ (n,D)=1}} g(n)  +  \sum\limits_j \frac{\overline{\chi_j}(a)}{\varphi(D)} \sum\limits_{n \le y} g(n)\chi_j(n)  \notag \\
& + O\Biggl( \frac{y}{\varphi(D)\log y} \prod\limits_{\substack{p \le D \\ (p,D)=1}} \left(1+\frac{|g(p)|}{p} \right) \left(\frac{\log y}{\log D}\right)^\alpha \Biggr) \notag
\end{align}
uniformly for $(a,D)=1$, $x^\gamma \le y \le x$.

In particular, the error term is $\ll D^{-1}y \left(\log D/\log y\right)^{1-\alpha}$.
\end{theorem} 


The accompanying paper, II, contains a complete proof of


\makeatletter{}

\begin{theorem}\label{JP_elliott_04_thm_02}
Let $3/2 \le Y \le x$.  Let $g$ be a multiplicative function that for positive constants $\beta, c, c_1$ satisfies $|g(p)| \le \beta$,
\begin{equation}
\sum\limits_{w < p \le x} \left(|g(p)| - c\right)p^{-1} \ge -c_1, \quad Y \le w \le x, \notag
\end{equation}
on the primes.  Suppose, further, that the series
\begin{equation}
\sum\limits_q |g(q)| q^{-1} (\log q)^\gamma, \quad \gamma = 1 + c\beta(c+\beta)^{-1}, \notag
\end{equation}
taken over the prime-powers $q = p^k$ with $k \ge 2$, converges.

Then with
\begin{equation}
\lambda = \min\limits_{|t| \le T} \sum\limits_{Y < p \le x} \left(|g(p)| - \Re g(p)p^{it}\right)p^{-1}, \notag
\end{equation}
\begin{equation}
\sum\limits_{n \le x} g(n) \ll x(\log x)^{-1} \prod\limits_{p \le x} \left(1+|g(p)|p^{-1}\right) \left( \exp\left(-\lambda c(c+\beta)^{-1} \right) + T^{-1/2} \right) \notag
\end{equation}
uniformly for $Y$, $x$, $T > 0$, the implied constant depending at most upon $\beta, c, c_1$ and a bound for the sum of the series over higher prime-powers.
\end{theorem} 


Although their arguments have features in common, the emphasis is sufficiently different that it seemed better to give each theorem its own presentation.  In Theorem \ref{JP_elliott_03_thm_01} functions are considered in packets; in Theorem \ref{JP_elliott_04_thm_02} they are considered singly.

Further, the present paper applies Theorem \ref{JP_elliott_04_thm_02} to illustrate a taxonomy of the characters appearing in Theorem \ref{JP_elliott_03_thm_01}.  The second paper collectively applies the various results to the study of automorphic forms.

To place the present paper in a wider context and to motivate the pervasive presence of the logarithmic function, we begin with an overview in the language of group representations.

$\mathbb{Q}^*$ will denote the multiplicative group of positive rationals.

The dual group of $\mathbb{Q}^*$, the direct product of denumerably many copies of $\mathbb{R}/\mathbb{Z}$, may be identified with the space of completely multiplicative functions $g$ with values in the complex unit circle, $|z|=1$.  We may topologize this space with a metric
\begin{equation}
\rho(g,h) = \left(\sum\limits_p p^{-\sigma} |g(p) - h(p)|^2 \right)^{1/2} \ge 0, \notag
\end{equation}
where $\sigma > 1$ and the sum is taken over the prime numbers.

It is sometimes convenient to employ a family of such metrics with differing values of $\sigma$.  Each metric is translation invariant and induces the standard topology attached to the dual of a locally compact abelian group.

Note that applied to general multiplicative functions with values in the complex unit disc, $|z| \le 1$, these are metrics on the equivalence classes of functions that coincide on the primes but not necessarily on the higher prime-powers.

Formally, or in an $L^2$ sense, a function $f:\mathbb{Q}^* \to \mathbb{C}$ gives rise to a Fourier transform
\begin{equation}
\hat{f}(g) = \int\limits_{\mathbb{Q}^*} f(r)g(r)\, d\eta(r) \notag
\end{equation}
where $\eta$ is purely atomic, assigning measure 1 to each positive rational, $r$.

Conversely, for each positive rational, $r$,
\begin{equation}
f(r) =  \int\limits_{\widehat{\mathbb{Q}^*}}\hat{f}(g)\overline{g(r)} \, dg \notag
\end{equation}
where $dg$ is the Haar measure on the compact dual group $\widehat{\mathbb{Q}^*}$, conveniently normalised to give measure 1 to the whole group.

In particular, restricting $r$ to the positive integers,
\begin{equation}
\sum\limits_{\substack{n \le x \\ n \equiv a \modd{D}}} f(n) = \int\limits_{\widehat{\mathbb{Q}^*}} \hat{f}(g)\sum\limits_{\substack{n \le x \\ n \equiv a \modd{D}}} \overline{g(n)} \, dg. \notag
\end{equation}

In the general study of arithmetic functions on residue classes, a specialisation to multiplicative functions with values in the complex unit disc then appears appropriate.

The one-dimensional unitary representations of the multiplicative group of positive reals into the invertible linear maps on $L^2(\mathbb{R})$ with respect to Lebesgue measure form a one-parameter group typically given by
\begin{equation}
t \mapsto (w \mapsto t^{i\alpha}w, w \in L^2(\mathbb{R})), \quad \alpha \in \mathbb{R}. \notag
\end{equation}
The infinitesimal generator of this group is given by the map $w \mapsto iw\log t$, $w \in L^2(\mathbb{R})$.

Similarly, we may view the representations
\begin{equation}
S_\tau: \mathbf{a} = (\dots, a_n, \dots) \mapsto (\dots, a_n n^{-i\tau}, \dots), \quad \tau \in \mathbb{R}, \notag
\end{equation}
as a one-parameter group of invertible maps into itself of the $\ell^2$ Hilbert space of complex-valued functions on the positive integers, norm $\left(\sum_{n=1}^\infty |a_n|^2 \right)^{1/2}$.  The Dirichlet series $\sum_{n=1}^\infty g(n)n^{-s}$, $s = \sigma + i\tau$, for a multiplicative function $g$ with values in the complex unit disc and $\sigma$ fixed at a positive value, may be identified with the orbit $S_\tau(\dots, g(n)n^{-\sigma}, \dots)$ whose elements generate an invariant subspace from which, as is shown in Lemma \ref{JP_elliott_03_lem_one_parameter_recovery}, $g$ may be recovered, up to translation, as any member with multiplicative coefficients.

The infinitesimal generator of the $S_\tau$ is given by
\begin{equation}
\mathbf{a} \mapsto -i(\dots, a_n\log n, \dots). \notag
\end{equation}
Viewed on the whole space $\ell^2$, this operator is not bounded.  However, in the present number-theoretical circumstances, the Euler product representation of the series $\sum_{n=1}^\infty g(n)n^{-s}$, whose existence is equivalent to the multiplicativity of the coefficient function $g$, affords a representation
\begin{align}
G(s) & = \exp\left(\sum\limits_p \log\left(1+g(p)p^{-s} + \cdots\right) \right) \notag \\
& = G(\sigma)\exp\left(i\tau G'(\sigma)/G(\sigma) + \cdots \right) \notag
\end{align}
corresponding to Stone's representation of $S_\tau$ in terms of its infinitesimal generator, rendering the ratio
\begin{equation}
G'(\sigma)/G(\sigma) = -\sum\limits_p \left(1+g(p)p^{-\sigma} + \cdots \right)^{-1} \sum\limits_{k=1}^\infty g(p^k)p^{-k\sigma} \log p^k, \notag
\end{equation}
corresponding to the action of the infinitesimal generator, manageable.

Note that the operator given by $G(\sigma) \mapsto G'(\sigma)/G(\sigma)$ only plays the role of the infinitesimal generator of the $S_\tau$, whose actual action is given by $G(\sigma) \mapsto -iG'(\sigma)$.

Although the expected logarithm function appears, the representing Dirichlet series is supported only on the prime-powers, which compensates.  In particular, corresponding to the factorisation $G(s)\left(G'(s)/G(s)\right)$ the function $g\log$ may be viewed as a convolution.  Otherwise expressed, we have factored $g\log$ in the group algebra of $\mathbb{Q}^*$.

Conveniently, for multiplicative functions with values in the complex unit disc
\begin{equation}
-G'(s)/G(s) = \sum\limits_p g(p)p^{-s}\log p - \psi'(s), \notag
\end{equation}
where
\begin{equation}
\psi(s) = \sum\limits_p \left(\log\left(1+g(p)p^{-s} + \cdots \right) - g(p)p^{-s} \right) \notag
\end{equation}
converges uniformly absolutely in each half-plane $\sigma \ge \sigma_1 > 1/2$, rendering the functions $\psi(s)$ and $\psi'(s)$ then uniformly bounded analytic.

Since the logarithmic function oscillates slowly on $\mathbb{R}$, a simple integration by parts may facilitate its removal.  As a consequence, the mean-value of $g$ may be attached to that of $g \log g$ and replaced by the specialisation of a bilinear form
\begin{equation}
x^{-1} \sum\limits_{\substack{p^km \le x \\ (p,m)=1}} g(m)g(p^k)\log p^k. \notag
\end{equation}

Extending each Dirichlet character to have the value 1 on primes that divide the corresponding modulus, we may regard the functions $g\chi_j$ appearing in Theorem \ref{JP_elliott_03_thm_01} as characters attached to the tensor product of one-dimensional representations of $\mathbb{Q}^*$ and $(\mathbb{Z}/D\mathbb{Z})^*$, respectively.  In a sense
\begin{equation}
\rho(1,g\chi_j) = \rho(\overline{\chi_j},g) \notag
\end{equation}
measures the distance of a typical tensor product from the trivial representation.

The wide spacing of Dirichlet characters on $\widehat{\mathbb{Q}^*}$ ensures that at most one can be near to a given general character, $g$.  Indeed, under mild constraints, and for a similar reason, at most one character attached to the tensor product of one-dimensional representations of $(\mathbb{Z}/D\mathbb{Z})^*$ and the group of positive reals can be near to a given multiplicative function, $g$, with values in the complex unit disc.  Theorem \ref{JP_elliott_03_thm_01} shows that, for practical purposes, estimation of the mean-value of $g$ on a given residue class involves only a small number of Dirichlet characters, their number independent of the size of the modulus.  In certain applications, such as to the study of primes in arithmetic progression, this uniformity is important.

Amongst other things, Lemma \ref{JP_elliott_03_lem_A1} asserts that the sets of primes on which various Dirichlet characters to reasonably sized moduli closely approximate a given function, $g$, are essentially disjoint, offering a tensor decomposition of the background representation of $\mathbb{Q}^*$ with details controlled by Dirichlet characters whose orders are also bounded independently of the associated moduli.

The sizes of $G(s,\chi) = \sum_{n=1}^\infty g(n)\chi(n)n^{-s}$ and the corresponding $\log G(s,\chi)$, Dirichlet $L$-functions attached to $g$ braided with a Dirichlet character $\chi$, are controlled by the Large Sieve, i.e. bounds for the spectra of appropriate self-adjoint operators acting upon spaces of functions supported on the integers or the primes, as the case may be.

From the viewpoint of group representations it is natural to consider functions in $L^2$ spaces.  The step from an $L^2$ estimate to an $L^\infty$ estimate that is presented in Theorem \ref{JP_elliott_03_thm_01} depends vitally upon attaching a bilinear form to the function $g$.

The maps $g \mapsto g(r)$, corresponding to the embedding of $\mathbb{Q}^*$ into its second dual, may be viewed as random variables with respect to the measurable sets of $\widehat{\mathbb{Q}^*}$, and the argument of the present paper appraised within the aesthetic of the theory of probability.  This brings into relief the need for maximal versions of the various inequalities arising.

The present argument takes place entirely in the half-plane of absolute convergence of the various Dirichlet series $\sum_{n=1}^\infty g(n)\chi(n)n^{-s}$; no analytic continuation is required of their sum functions $G(s,\chi)$; boundary value behaviour is essentially classified.  This affords applications to problems that are otherwise currently out of reach.

Further remarks, including those of an historical nature, may be found in the concluding section of this paper.


\section{Inequalities of Large Sieve type}


Besides error terms, the following estimates of operator norms will control the size of $L$-functions and their logarithms.


\makeatletter{}

\begin{lemma}\label{JP_elliott_03_lem_I1}
Let $0 < \varepsilon < 1$.  The inequality
\begin{equation}
\sum\limits_{j=1}^J \max\limits_{v-u \le H} \left|\sum\limits_{\substack{u < n \le v \\ (n,Q)=1}} a_n\chi_j(n)\right|^2 \ll \left(H\prod\limits_{\substack{p \mid Q \\ p \le H}} \left(1-\frac{1}{p}\right) + JH^\varepsilon D^{1/2}\log D\right)\sum\limits_{n=1}^\infty |a_n|^2 \notag
\end{equation}
where the $\chi_j$ are distinct Dirichlet characters $\modd{D}$, $D \ge 2$, $Q$ a positive integer, $H \ge 0$, holds for all square-summable complex numbers $a_n$, the implied constant depending at most upon $\varepsilon$.
\end{lemma} 

There are several ways to establish this result, which is of Maximal Gap Large Sieve type.  An application of Cauchy's inequality shows that we may assume $D^{1/2}$ not to exceed $H$.  We may also include in $Q$ the prime-divisors of $D$.  Note that $\sum\limits_{\substack{p \mid Q \\ x < p \le x^2}} p^{-1} \le \sum\limits_{x < p \le x^2} p^{-1} \ll 1$ uniformly in $x \ge 1$.

\makeatletter{}

\noindent \emph{Proof of Lemma \ref{JP_elliott_03_lem_I1}}.  With $0 \le v_j-u_j \le H$, define
\begin{equation}
t_j(n) = \begin{cases} \chi_j(n) & \text{if} \ u_j < n \le v_j, \\ 0 & \text{otherwise}, \end{cases} \qquad j = 1, \dots, J. \notag
\end{equation}
For any real $\lambda_d$, $d\mid Q$, constrained by $\lambda_1 = 1$, the dual form
\begin{equation}
S = \sum\limits_{\substack{n \le x \\ (n,Q)=1}} \left|\sum\limits_{j=1}^J c_jt_j(n)\right|^2 \notag
\end{equation}
does not exceed
\begin{equation}
\sum\limits_{n \le x} \left(\sum\limits_{d \mid (n,Q)} \lambda_d \right)^2 \left|\sum\limits_{j=1}^J c_jt_j(n)\right|^2 \notag
\end{equation}
\begin{equation}
= \sum\limits_{d_i \mid Q} \lambda_{d_1}\lambda_{d_2} \sum\limits_{j,k=1}^J c_j \overline{c_k} \sum\limits_{ n \equiv 0 \modd{[d_1,d_2]}} t_j(n) \overline{t_k(n)}, \notag
\end{equation}
where $[d_1,d_2]$ denotes the least common multiple of $d_1$ and $d_2$.  For those terms with $j \ne k$, the innermost sum has the form
\begin{equation}
\chi_j\overline{\chi_k}\left([d_1,d_2]\right)\sum\limits_{m} \chi_j\overline{\chi_k}(m) \notag
\end{equation}
with the integers $m$ over an interval and is, by a classical result of P\'olya and Vinogradov, $O\left(D^{1/2}\log D\right).$  The corresponding contribution to $S$ is
\begin{equation}
\ll JD^{1/2}\log D \left(\sum\limits_{d \mid Q} |\lambda_d|\right)^2 \sum\limits_{j=1}^J |c_j|^2. \notag
\end{equation}
For those terms with $j=k$ we reform the square in the $\lambda_d$ to gain a contribution
\begin{equation}
\sum\limits_{j=1}^J |c_j|^2 \sum\limits_{n \le x} \left(\sum\limits_{d\mid (n,Q)} \lambda_d\right)^2 \left|t_j(n)\right|^2. \notag
\end{equation}
Since $\left|t_j(n)\right| \le 1$, the innersum over $n$ does not exceed
\begin{equation}
\sum\limits_{d_i \mid Q} \lambda_{d_1}\lambda_{d_2} \sum\limits_{\substack{u_j < n \le v_j+H \\ n \equiv 0 \modd{[d_1,d_2]}}} 1 = H\sum\limits_{d_i \mid Q} \lambda_{d_1}\lambda_{d_2}[d_1,d_2]^{-1} + O\left(\left(\sum\limits_{d\mid Q} |\lambda_d|\right)^2\right). \notag
\end{equation}
We may follow the standard appeal to the method of Selberg, c.f. Elliott \cite{Elliott1979}, Chapter 2, with $\lambda_d=0$ if $d>H^{\varepsilon/2}$ which, in particular, gives $|\lambda_d| \le 1$ for all remaining $\lambda_d$.

As a consequence,
\begin{equation}
S \ll \left(H\prod\limits_{\substack{p \mid Q \\ p \le H}} \left(1-\frac{1}{p}\right) + JH^\varepsilon D^{1/2}\log D\right)\sum\limits_{j=1}^J |c_j|^2. \notag
\end{equation}
Dualising gives the inequality of Lemma \ref{JP_elliott_03_lem_I1}. 

Specialising $H$ to $x$, $Q$ to $\prod_{p \le x^\varepsilon} p$ we obtain the

\begin{corollary} Let $0 < \varepsilon < 1$.  The inequality
\begin{equation}
\sum\limits_{\chi} \left| \sum\limits_{q \le x} a_q \chi(q)\right|^2 \ll \left(\frac{x}{\log x} + x^\varepsilon D^{3/2}\log D \right) \sum\limits_{q \le x} |a_q|^2 \notag
\end{equation}
where $\chi$ traverses the characters $\chi\modd{D}$, and $q$ the prime-powers, holds for all complex numbers $a_q$, real $x \ge 2$.
\end{corollary}


\makeatletter{}

\begin{lemma} \label{JP_elliott_03_lem_I2}
There is a real $c$ such that
\begin{equation}
\sum\limits_{j=1}^k \max\limits_{D \le w \le y \le x} \max\limits_{\sigma \ge 1, |t| \le T} \left| \sum\limits_{w < p \le y} a_p \chi_j(p) p^{-s} \right|^2 \le 4(L + k\Delta) \sum\limits_{D < p \le x} |a_p|^2 p^{-1}, \notag
\end{equation}
with $s = \sigma + it$, $\sigma = \Re(s)$, $L = \sum_{D< p \le x} p^{-1}$, $\Delta = \log(\log T/\log D) + c$, uniformly for $a_p$ in $\mathbb{C}$ and distinct characters $\chi_j \modd{D}$, $j = 1, \dots, k$, $x \ge D \ge 2$, $T \ge D$.
\end{lemma} 

Whilst a version of Lemma \ref{JP_elliott_03_lem_I2} may be deduced from the Corollary to Lemma \ref{JP_elliott_03_lem_I1} using an integration by parts, the dependence of the resulting bound upon the size of $|t|$ is severe.  The amelioration supplied by the following result is vital.


\makeatletter{}

\begin{lemma} \label{JP_elliott_03_lem_I3}
With a certain constant $c$,
\begin{equation}
\Re \sum\limits_{y < p \le w} \chi(p) p^{-s} \le \log(\log T/\log D)+c \notag
\end{equation}
uniformly for nonprincipal characters $\chi \modd{D}$, $\sigma \ge 1$, $|t| \le T$, $w \ge y \ge D$, $T \ge D \ge 2$.
\end{lemma} 

\makeatletter{}

\noindent \emph{Proof of Lemma \ref{JP_elliott_03_lem_I3}}.  We shall appeal to the following bound for Dirichlet $L$-functions:
\begin{equation}
\left| \frac{L(\sigma_1 + it, \chi)}{L(\sigma_2 + it, \chi)} \right| \le \left(D(|t|+2) \right)^{c_1(\sigma_2 - \sigma_1)}, \notag
\end{equation}
valid for $1 < \sigma_1 \le \sigma_2$ with an absolute constant $c_1$.  Two proofs via analytic functions may be found in \cite{elliottMFoAP6}, the restriction to $\sigma_2 \le 2$ there unnecessary;  an alternative elementary argument via a sieve is given in \cite{elliott2002millenium}.

For $\beta > 1$, $p^{-1} - p^{-\beta} \le (\beta-1)p^{-1} \log p$ holds and integration by parts together with the well-known Chebyshev bound $\pi(x) \ll x/\log x$ shows the sums
\begin{equation}
\sum\limits_{p > \exp\left((\beta-1)^{-1}\right)} p^{-\beta}, \quad \sum\limits_{p \le \exp\left((\beta-1)^{-1}\right)} \left(p^{-1} - p^{-\beta} \right), \quad (\beta-1) \sum\limits_{p \le \exp\left((\beta-1)^{-1}\right)} p^{-1}\log p \notag
\end{equation}
to be uniformly bounded.

With $\beta = 1+ (\log w)^{-1}$,
\begin{align}
\sum\limits_{p \le w} \chi(p)p^{-s}&  - \sum\limits_p \chi(p) p^{-(\beta + s-1)} \notag \\
&  = \sum\limits_{p \le \exp\left((\beta-1)^{-1}\right)} \chi(p)p^{-s+1}\left(p^{-1} - p^{-\beta}\right) - \sum\limits_{p > \exp\left((\beta-1)^{-1} \right)} \chi(p)p^{-(\beta+s-1)} \notag
\end{align}
is bounded uniformly for $w \ge 2$; likewise when $w$ is replaced by $y$.

From the Euler product representation
\begin{equation}
L(z,\chi) = \prod\limits_p \left(1-\chi(p)p^{-z} \right)^{-1} = \exp\left(\sum\limits_p \chi(p)p^{-z} + O(1) \right), \notag
\end{equation}
valid for $\Re(z) > 1$,
\begin{equation}
\exp\left(\Re \sum\limits_{y < p \le w} \chi(p)p^{-s} \right) \ll \left|\frac{ L\left(s + \frac{1}{\log w}, \chi\right)}{L\left(s + \frac{1}{\log y}, \chi\right)} \right| \ll \left(D(|t|+2) \right)^{c_1/\log y}. \notag
\end{equation}
Bearing in mind the restriction $y \ge D$, if $y \ge T$ further holds then this upper bound does not exceed an absolute constant.  Otherwise we apply the same argument to the range $T < p \le w$ and note that
\begin{equation}
\left| \Re \sum\limits_{y < p \le T} \chi(p)p^{-s} \right| \le \sum\limits_{D < p \le T} p^{-1} \le \log(\log T/\log D) + c_2. \notag
\end{equation}

Taking logarithms completes the proof.

\makeatletter{}

\noindent \emph{Proof of Lemma \ref{JP_elliott_03_lem_I2}}.  Since the sum $\sum_{D<p \le x} |a_p|p^{-\sigma}$ approaches zero as $\sigma \to \infty$, the innermost maximum may be taken over a bounded rectangle.

For reals $y_j$, $w_j$, $D \le y_j \le w_j$, $\sigma_j \ge 1/2$, $t_j$, $|t_j| \le T$, define
\begin{equation}
\delta_{j,p} = \begin{cases} \chi_j(p)p^{-\sigma_j - it_j} & \text{if} \ y_j <p \le w_j, \\ 0 & \text{otherwise,} \end{cases} \notag
\end{equation}
$j = 1, \dots, k$, and consider the inequality
\begin{equation}
\sum\limits_{D<p \le x} \left|\sum\limits_{j=1}^k b_j \delta_{j,p} \right|^2 \le \lambda \sum\limits_{j=1}^k |b_j|^2, \notag
\end{equation}
where the $b_j$ are for the moment real and nonnegative.  The expanded sum is
\begin{equation}
\sum\limits_{j=1}^k b_j^2 \sum\limits_{y_j < p \le w_j} p^{-2\sigma_j} + 2\sum\limits_{1 \le j < \ell \le k} b_jb_\ell \,\Re\sum\limits_{D < p \le x} \chi_j \overline{\chi_\ell}(p)p^{-\sigma_j - \sigma_\ell-it_j+it_\ell}. \notag
\end{equation}
An appeal to Lemma \ref{JP_elliott_03_lem_I3} followed by an application of the Cauchy-Schwarz inequality shows that with a suitable choice for $c$ we may take $\lambda = L + k\Delta$.

If now $b_j$ is complex, we represent it as a sum
\begin{equation}
\max(\text{Re} \ b_j,0) + \min(\text{Re} \ b_j,0) + i\max(\text{Im} \ b_j,0)+i\min(\text{Im} \ b_j,0) \notag
\end{equation}
and correspondingly partition the innersum over $j$.  Since the coefficients in each subsum all have the same argument, a second application of the Cauchy-Schwarz inequality allows us to conclude that with $\lambda = 4(L + k\Delta)$ the above inequality holds for all complex $b_j$.

Dualising:
\begin{equation}
\sum\limits_{j=1}^k \left|\sum\limits_{D<p \le x} a_p \delta_{j,p} \right|^2 \le 4(L + k\Delta)\sum\limits_{D<p \le x} |a_p|^2 \notag
\end{equation}
for all complex $a_p$.

Replacing $a_p$ by $a_pp^{-\frac{1}{2}}$ completes the proof. 


\section{Exceptional Characters}


We may appreciate Lemma \ref{JP_elliott_03_lem_I2} by applying it to the Dirichlet series
\begin{equation}
G_{I,j}(s) = \sum\limits_{n=1}^\infty g(n)\chi_j(n)n^{-s} \notag
\end{equation}
where the multiplicative function $g$, with values in the complex unit disc, vanishes on the primes outside the interval $I$, and $\chi_j$ is a nonprincipal character $\modd{D}$.

For $\sigma > 1$, $I \subseteq (D,x]$, in terms of the principal value of the logarithm,
\begin{equation}
\log G_{I,j}(s) = \sum\limits_{p \in I} g(p)\chi_j(p) p^{-s} + O\left(\sum\limits_{p > D} p^{-2} \right). \notag
\end{equation}
If $H$ is the semi-strip $\Re(s) > 1$, $|\Im(s)| \le T$, then
\begin{equation}
\sum\limits_{j=1}^k \max\limits_{I \subseteq (D,x]} \max\limits_{s \in H} \left| \log G_{I,j}(s) \right|^2 \le 8(L + k\Delta)L + O\left((\log D)^{-1} \right). \notag
\end{equation}
%

Given $\alpha>0$, $T \le \exp\left(\log D(\log x/\log D)^{\alpha^2/9} \right)$, we call a nonprincipal character $\chi_j \modd{D}$ for which
\begin{equation}
\max\limits_{I \subseteq (D,x]} \max\limits_{s \in H} \left| \log G_{I,j}(s) \right| < \alpha \log (\log x/\log D) \notag
\end{equation}
fails, \emph{exceptional relative to the triple} $( \alpha, D, x )$ or, more shortly, \emph{exceptional}.

Thus, given $\alpha > 0$, with the exception of $O\left(\alpha^{-2}\right)$ characters,
\begin{equation}
\left( \frac{\log x}{\log D} \right)^{-\alpha} < \left| G_{I,j}(s) \right| < \left(\frac{\log x}{\log D}\right)^\alpha \notag
\end{equation}
uniformly for $s$ in $H$, $I$ in $(D,x]$.

Then either there are $\ll \alpha^{-2}$ such characters attached to the modulus $D$, or $D$ exceeds a certain fixed power of $x$, the power depending upon $\alpha$ only.

\emph{Remarks}.  If $0 < \gamma < 1$, $D \le N^\gamma \le x \le N$, then a character $\chi_j$ that is nonexceptional relative to $( \alpha, D, N)$ is essentially nonexceptional relative to $(\alpha, D, x)$, since each subset of $(D,x]$ is a subset of $(D,N]$ and
\begin{equation}
\left|G_{I,j}(s) \right| < \left(\frac{\log N}{\log D} \right)^\alpha \le \left( \frac{\log x}{\gamma \log D} \right)^\alpha, \notag
\end{equation}
with an analogous lower bound.

Note that within similar tolerances the same result applies to $g\mu$, obtained by braiding $g$ with the M\"obius function.


\section{First Waystation:  an $L^2$ Theorem \ref{JP_elliott_03_thm_01}} \label{JP_elliott_03_sec_first_waystation}


\makeatletter{}

\begin{lemma} \label{JP_elliott_03_lem_I4}
Let $\gamma < 1$, $\delta < 1$, $c>3/2$, $N$ be positive real numbers, $D$ an integer, $2 \le D \le N^{1/c}$.

Given a completely multiplicative function $g$ that vanishes on the primes in $(1,D^c]$, there is a set $\mathcal{J}$ of Dirichlet characters $\modd{D}$, of cardinality bounded in terms of $\delta$ alone, so that
\begin{equation}
\sum\limits_{\chi \notin \mathcal{J}} \max\limits_{2 \le y \le t} \left| \sum\limits_{n \le y} g(n)\chi(n) \right|^2 \ll \left(\frac{t}{\log t}\right)^2 \left(\frac{\log x}{\log D}\right)^\delta \notag
\end{equation}
uniformly for $D^c \le t \le x$, $N^\gamma \le x \le N$.
\end{lemma} 

We may interpret this result by means of

\makeatletter{}

\begin{lemma} \label{JP_elliott_03_lem_I5}
Let $\chi_j$, $j \in \mathcal{J}$, be a collection of Dirichlet characters $\modd{D}$.  Let
\begin{equation}
\mathcal{L}(a) = \sum\limits_{n \equiv a \modd{D}} b_n - \sum\limits_{j \in \mathcal{J}} \frac{\overline{\chi_j}(a)}{\varphi(D)} B_j \notag
\end{equation}
where at most finitely many of the complex numbers $b_n$ are nonzero and
\begin{equation}
B_j = \sum\limits_n b_n\chi_j(n). \notag
\end{equation}

Then
\begin{equation}
\varphi(D)\sum\limits_{\substack{a=1 \\ (a,D)=1}}^D \left| \mathcal{L}(a) \right|^2 = \sum\limits_{j \notin \mathcal{J}} |B_j|^2. \notag
\end{equation}
\end{lemma} 

\makeatletter{}

\noindent \emph{Proof of Lemma \ref{JP_elliott_03_lem_I5}}.  From the orthogonality of characters
\begin{equation}
\mathcal{L}(a) = \sum\limits_{j \notin \mathcal{J}} \frac{\overline{\chi_j}(a)}{\varphi(D)} B_j. \notag
\end{equation}
Hence
\begin{align}
\sum\limits_{\substack{a=1 \\ (a,D)=1}}^D |\mathcal{L}(a)|^2 & = \sum\limits_{\substack{a=1 \\ (a,D)=1}}^D \frac{1}{\varphi(D)^2} \sum\limits_{j_1, j_2 \notin \mathcal{J}} \overline{\chi_{j_1}}\chi_{j_2}(a)B_{j_1}\overline{B_{j_2}} \notag \\
& = \frac{1}{\varphi(D)^2} \sum\limits_{j_1, j_2 \notin \mathcal{J}} \sum\limits_{\substack{a=1 \\ (a,D)=1}}^D \overline{\chi_{j_1}}\chi_{j_2}(a) \notag \\
& = \frac{1}{\varphi(D)} \sum\limits_{j \notin \mathcal{J}} |B_j|^2, \notag
\end{align}
the final step by a further appeal to the orthogonality of characters. 

For an arithmetic function $f$, define
\begin{equation}
Y(f,a,x) = \sum\limits_{\substack{n \le x \\ n \equiv a \modd{D}}} f(n) - \sum\limits_{j \in \mathcal{J}} \frac{\overline{\chi_j}(a)}{\varphi(D)} \sum\limits_{n \le x} f(n)\chi_j(n). \notag
\end{equation}
Then Lemma \ref{JP_elliott_03_lem_I4} in particular asserts that
\begin{equation}
\varphi(D) \sum\limits_{\substack{a=1 \\ (a,D)=1}}^D |Y(g,a,t)|^2 \ll \left(\frac{t}{\log t}\right)^2 \left(\frac{\log x}{\log D}\right)^\delta \notag
\end{equation}
uniformly for $D^c \le t \le x$, $N^\gamma \le x \le N$, a version of Theorem \ref{JP_elliott_03_thm_01} that belongs formally to the $L^2$ space of functions on the group of reduced residues $\modd{D}$.

For a complex-valued arithmetic function $f$ and real $x$ define
\begin{equation}
N(x) = N(x,f) = \sum\limits_{n \le x} f(n)\log n, \qquad M(x) = M(x,f) = \sum\limits_{n \le x} f(n). \notag
\end{equation}
When it is clear from the context, explicit reference to the function $f$ may be omitted.

$\Lambda(n)$ denotes von Mangoldt's function, $\log p$ when $n$ is a power of a prime $p$, zero otherwise.

\makeatletter{}

\begin{lemma} \label{JP_elliott_03_lem_I6}
Let $r>0$, $y = w-w(\log w)^{-r}$.  Then
\begin{equation}
N(w,g) \ll w\int_2^w \frac{|N(u)|}{u^2\log u}\, du + \sum\limits_{d \le (\log w)^{2r}} \frac{d\Lambda(d)|g(d)|}{w-y} \int_{y/d}^{w/d} |M(u)| \, du + E_0, \notag
\end{equation}
where
\begin{equation}
E_0 = E_0(g) = \max\limits_{y \le t \le w} \left| \sum\limits_{y < n \le t} g(n)\log n \right| \notag
\end{equation}
uniformly for $w \ge 2$, for all completely (or exponentially) multiplicative functions $g$ with values in the complex unit disc.
\end{lemma} 

\noindent \emph{Proof of Lemma \ref{JP_elliott_03_lem_I6}}.  This result may be found as Lemma 1 of \cite{elliottMFoAP4}.

The following argument, which formally applies to any complex-valued arithmetic function $f$, will be implicitly employed several times.

Integrating by parts, $F(s)$, $s = \sigma + i\tau$, $\sigma = \Re(s)>0$, the formal sum function of the Dirichlet series $\sum_{n=1}^\infty f(n)n^{-s}$, has representations
\begin{equation}
s^{-1}F(s) = \int_1^\infty y^{-s-1}M(y,f)\, dy = \int_1^\infty M(e^w)e^{-w\sigma}e^{-iw\tau}\, dw.  \notag
\end{equation}
Viewed as functions of $\tau$ and $w$ respectively, $\left(s\sqrt{2\pi}\right)^{-1}F(s)$ and $M(e^w)e^{-w\sigma}$ are Fourier transforms.  By Plancherel's theorem:

\makeatletter{}

\begin{lemma} \label{JP_elliott_03_lem_I7}
\begin{equation}
\int_{-\infty}^\infty \left| \frac{F(s)}{s}\right|^2 \,d\tau = 2\pi \int_1^\infty \frac{|M(y)|}{y^{2\sigma+1}} \, dy \notag
\end{equation}
provided one of the integrals exists in an $L^2$ sense.
\end{lemma} 

\makeatletter{}

\noindent \emph{Proof of Lemma \ref{JP_elliott_03_lem_I4}}.  For simplicity of presentation we shall establish the result with $3\delta$ in place of $\delta$.  Temporarily denoting $\log t$ by $\ell$, it follows from Lemma \ref{JP_elliott_03_lem_I6} that for any set of distinct characters $\chi_j \modd{D}$,
\begin{equation}
S(t) = \sum\limits_j \max\limits_{2 \le w \le t} \left|N(w,g\chi_j)\right|^2 \ll \sum\limits_{k=1}^3 F_k, \notag
\end{equation}
where
\begin{align}
F_1 & = \sum\limits_j \left(t \int_2^t \frac{|N(u)|}{u^2\log u} \, du \right)^2, \notag \\
F_2 & = \sum\limits_j \left( \sum\limits_{d \le \ell^{2r}} \Lambda(d)|g(d)| \max\limits_{u \le d^{-1}t} \left|M(u,g\chi_j)\right| \right)^2, \notag \\
F_3 & = \sum\limits_j \left( \max\limits_{\substack{ v-u \le \ell^{-r}t \\ v \le t}} \left| \sum\limits_{u < n \le v} g(n)\chi_j(n)\log n \right| \right)^2. \notag
\end{align}

We consider these expressions in reverse order.  Assuming $c>3/2$, an application of Lemma \ref{JP_elliott_03_lem_I1} with $r=2$ and $\varepsilon$ sufficiently small shows that
\begin{equation}
F_3 \ll \left(t(\log t)^{-2} + D^{3/2}\log D \log t \right)\sum\limits_{n \le t} (\log n)^2 \ll t^2. \notag
\end{equation}

With $\kappa = \sum_{D < d \le \ell^{2r}} d^{-1}\Lambda(d)$, after the Cauchy-Schwarz inequality a similar application of Lemma \ref{JP_elliott_03_lem_I1} delivers
\begin{align}
F_2 & \ll \kappa \sum\limits_{D^c < d \le \ell^{2r}} d\Lambda(d) \sum\limits_j \max\limits_{u \le d^{-1}t} |M(u,g\chi_j)|^2 \notag \\
& \ll \kappa \sum\limits_{D^c < d \le \ell^{2r}} d\Lambda(d) \left(d^{-1}t + (d^{-1}t)^\varepsilon D^{3/2}\log D \right) d^{-1}t \notag \\
& \ll (\kappa t)^2 \ll \left(t\log(\log t/\log D)\right)^2. \notag
\end{align}
Note that if $D^c \ge (\log t)^4$, then $F_2 = 0$.

Towards $F_1$ we also note that since $g$ vanishes on the integers not exceeding $D^c$, appeal to Lemma \ref{JP_elliott_03_lem_I1} delivers the bound
\begin{align}
\sum\limits_j |N(u,g\chi_j)|^2 & \ll \left(u(\log D)^{-1} + u^\varepsilon D^{3/2}\log D\right) \sum\limits_{n \le u} |g(n)\log n|^2  \notag \\
& \ll \left(u \log u/\log D\right)^2, \notag
\end{align}
uniformly for $u \ge 2$.

Let $\theta$ be a real number, $0 < \theta < 1$; its value will ultimately depend at most upon $c$, $\delta$, $\gamma$.  Let
\begin{equation}
S_1(t) = t^2 \sum\limits_j \left( \int_{t^\theta}^t \frac{|N(u,g\chi_j)|}{u^2 \log u} \, du \right)^2 \notag
\end{equation}
and $S_2(t)$ the similar expression with range of integration $D^c < u \le t^\theta$.

We begin with the Cauchy-Schwarz inequality:
\begin{equation}
S_1(t) \le t^2 \int_{t^\theta}^t \frac{du}{u(\log u)^2} \sum\limits_j \int_{t^\theta}^t \frac{|N(u,g\chi_j)|^2}{u^3} \, du. \notag
\end{equation}
Setting $\sigma = 1+(\log t)^{-1}$, $u^{-3} \le e^2 u^{-2\sigma-1}$ holds over the range $1 \le u \le t$, hence
\begin{align}
S_1(t) & \le \frac{(et)^2}{\theta \log t} \sum\limits_j \int_{t^\theta}^t \frac{|N(u,g\chi_j)|^2}{u^{2\sigma+1}} \, du \notag \\
& \le \frac{(et)^2}{2\pi\theta \log t} \sum\limits_j \int_{-\infty}^\infty \left|\frac{G'(s,\chi_j)}{s} \right|^2 \, d\tau, \notag
\end{align}
the second step by appeal to Lemma \ref{JP_elliott_03_lem_I7}.

For $T > 0$, to be chosen shortly, let
\begin{equation}
L_1 = \frac{t^2}{\theta \log t} \sum\limits_j \int_{|\tau| \le T} \left| \frac{G'(s,\chi_j)}{s}\right|^2 d\tau \notag
\end{equation}
and $L_2$ the similar expression with range of integration $|\tau| > T$.

A further application of Lemma \ref{JP_elliott_03_lem_I7} shows that
\begin{equation}
\sum\limits_j \int_{-\infty} ^\infty \left|\frac{G'(s,\chi_j)}{s}\right|^2 \, d\tau \ll \int_2^\infty \sum\limits_j |N(u,g\chi_j)|^2 \frac{du}{u^{2\sigma+1}}, \notag
\end{equation}
which our initial remark guarantees to be
\begin{equation}
\ll \int_2^\infty \left(\frac{\log u}{\log D}\right)^2 \frac{du}{u^{2\sigma-1}} \ll \frac{(\log t)^3}{(\log D)^2}. \notag
\end{equation}
In particular,
\begin{equation}
\sum\limits_j \int_{|\tau| \le 1} |G'(s,\chi_j)|^2 \, d\tau \ll (\log t)^3(\log D)^{-2}. \notag
\end{equation}

Since we can replace $g(n)$ by $g(n)n^{-i\lambda}$ for any real $\lambda$ without affecting the hypothesis that $g$ vanishes on the primes up to $D^c$, this last inequality holds when the integration is over any interval $|\tau - m| \le 1$, $m$ an integer.

As a consequence,
\begin{equation}
L_2 \ll \frac{t^2}{\theta \log t} \sum\limits_{|m| > T} \frac{1}{m^2} \frac{(\log t)^3}{(\log D)^2} \ll \frac{1}{\theta T} \left(\frac{t\log t}{\log D}\right)^2. \notag
\end{equation}

\emph{We now remove from consideration those characters $\chi_j$ that are exceptional relative to the triple $(\delta, D, N)$}.

With the implicit constraint that $T$ not exceed $\exp\left(\log D(\log x/\log D)^{\delta^2/9}\right)$, the factorisation $G'(s,\chi_j) = G(s,\chi)\left(G'(s,\chi)/G(s,\chi)\right)$ enables us to assert that
\begin{equation}
L_1 \ll \frac{t^2}{\theta\log t} \left(\frac{\log x}{\log D}\right)^{2\delta} \sum\limits_j \int_{-\infty}^\infty \left|\frac{G'(s,\chi_j)}{sG(s,\chi_j)}\right|^2 \, d\tau. \notag
\end{equation}

Appeal to the representation
\begin{equation}
- \frac{G'(s,\chi)}{G(s,\chi)} = \sum\limits_{n=1}^\infty \frac{g(n)\chi(n)\Lambda(n)}{n^s} \notag
\end{equation}
and application of Lemma \ref{JP_elliott_03_lem_I7} provides the sum over $j$ with the alternative representation
\begin{equation}
2\pi \int_1^\infty \sum\limits_j \left| \sum\limits_{n \le y} g(n)\chi_j(n) \Lambda(n) \right|^2 y^{-2\sigma-1} \, dy. \notag
\end{equation}
Here the integrand is zero unless $y > D^c$, when the corollary to Lemma \ref{JP_elliott_03_lem_I1} shows it to be
\begin{equation}
\ll \left(y(\log y)^{-1} + y^\varepsilon D^{3/2}\log D \right) \sum\limits_{n \le y} \Lambda(n)^2 y^{-2\sigma-1} \ll y^{-2\sigma + 1}. \notag
\end{equation}
Consequently $L_1 \ll t^2\theta^{-1}(\log x/\log D)^{2\delta}$ and
\begin{equation}
S_1(t) \ll t^2\theta^{-1} (\log x/\log D)^{2\delta} + t^2(\theta T)^{-1}(\log t/\log D)^2. \notag
\end{equation}
Choosing $T = c_0(\log N/\log D)^4$ with a suitably small constant $c_0$, we may omit the second of the bounding terms in favour of the first.

The sum $S_2(t)$ is treated indirectly.  For $D^c \le u \le x$ define
\begin{equation}
H(u) = \max\limits_{D^c \le w \le u} \frac{1}{w^2(\log w)^\delta} \sum\limits_j \max\limits_{2 \le y \le w} |N(y,g\chi_j)|^2. \notag
\end{equation}

After an application of the Cauchy-Schwarz inequality
\begin{align}
S_2(t) & \le t^2 \int_{D^c}^{t^\theta} \frac{du}{u(\log u)^{1-\delta/2}} \int_{D^c}^{t^\theta} \frac{H(u)}{u(\log u)^{1-\delta/2}} \, du \notag \\
& \ll t^2(\theta\log t)^\delta H(t^\theta) \ll t^2(\theta\log t)^\delta H(x^\theta), \notag
\end{align}
uniformly for $D^c \le t \le x$.  Then uniformly for $D^c \le t \le x$, $N^\gamma \le x \le N$,
\begin{equation}
S(t) \ll t^2(\theta\log t)^\delta H(x^\theta) + t^2\theta^{-1}(\log x/\log D)^{2\delta}, \notag
\end{equation}
the contribution from terms $F_2$, $F_3$ having been absorbed.  In particular,
\begin{equation}
H(x) \ll \theta^\delta H(x^\theta) + \theta^{-1} (\log x/\log D)^{2\delta} (\log D)^{-\delta}. \notag
\end{equation}
Fixing $\theta$ at a sufficiently small value, independent of $x$, $N$, $D$, we may omit the term involving $H(x^\theta)$ and arrive at
\begin{equation}
\sum_j \max\limits_{2 \le y \le t} |N(y,g\chi_j)|^2 \ll t^2(\log x/\log D)^{3\delta}, \notag
\end{equation}
with the same uniformities in $t$ and $x$.

To strip the logarithm in $N(y)$ we apply Lemma \ref{JP_elliott_03_lem_I1}:
\begin{align}
\sum\limits_j \max\limits_{2 \le y \le t}  \left| \sum\limits_{n \le y} g(n)\chi_j(n)\log(t/n) \right|^2 
 \ll \left(t+ t^\varepsilon D^{3/2}\log D \right) \sum\limits_{n \le t}  \log (t/n)^2 \ll t^2, \notag
\end{align}
the Cauchy-Schwarz inequality; and divide by $(\log t)^2$.

This completes the proof of Lemma \ref{JP_elliott_03_lem_I4}.


\section{Second Waystation:  from $L^2$ to $L^\infty$}


\makeatletter{}

\begin{lemma} \label{JP_elliott_03_lem_I8}
Let $c>3/2$.  Theorem \ref{JP_elliott_03_thm_01} is valid for completely multiplicative functions that vanish on the primes up to $D^c$.
\end{lemma} 

\makeatletter{}

\noindent \emph{Proof of Lemma \ref{JP_elliott_03_lem_I8}}.  We again introduce a logarithm, then remove it.

The representation $\log n = \sum_{d \mid n} \Lambda(d)$  enables the convolution factorisation $g\log = g*g\Lambda$, hence a representation
\begin{equation}
Y(g\log, a, t) = \sum\limits_{d \le t} g(d)\Lambda(d) Y(g,a\overline{d}, td^{-1}), \notag
\end{equation}
where $d \overline{d} \equiv 1 \modd{D}$.

The contribution from the terms with $tD^{-2c} < d \le t$ is
\begin{align}
& \ll \sum\limits_{tD^{-2c} < d \le t} \Lambda(d) \left( \sum\limits_{\substack{ m \le td^{-1} \\ m \equiv a\overline{d} \modd{D}}} g(m) + \frac{1}{\varphi(D)} \sum\limits_{m \le td^{-1}} |g(m)| \right) \notag \\
& \ll \sum\limits_{m \le D^{2c}} |g(m)| \left( \sum\limits_{\substack{d \le tm^{-1} \\ d \equiv a\overline{m} \modd{D}}} \Lambda(d) + \frac{1}{\varphi(D)} \sum\limits_{d \le tm^{-1}} \Lambda(d) \right) \notag
\end{align}
and, by the Brun-Titchmarsh theorem,
\begin{equation}
\ll \frac{t}{\varphi(D)} \sum\limits_{m \le D^{2c}} \frac{|g(m)|}{m} \ll \frac{t}{\varphi(D)} \prod\limits_{D^c < p \le D^{2c}} \left(1+ \frac{1}{p}\right) \ll \frac{t}{\varphi(D)}. \notag
\end{equation}

Let $0 < \beta < 1$.  We cover the interval $(D^c, tD^{-2c}]$ by adjoining intervals $(U,2U]$ and each such interval by adjoining subintervals $(V,V+U^\beta]$.  The remaining terms in the sum representing $Y(g\log, a, t)$ contribute
\begin{equation}
\ll \sum\limits_U \sum\limits_V \sum\limits_{V < d \le V+U^\beta} |g(d)| \Lambda(d) \left| Y(g,a\overline{d}, td^{-1}) \right|. \notag
\end{equation}

Replacing $td^{-1}$ in the innermost sum by $tV^{-1}$ introduces an error of
\begin{align}
& \ll \sum\limits_U \sum\limits_V \sum\limits_{V < d \le V+U^\beta} \Lambda(d) \left( \frac{1}{\varphi(D)} \left( \frac{t}{V} - \frac{t}{d} \right)+1\right) \notag \\
& \ll \sum\limits_{D^c < d \le 2tD^{-2c}} \frac{t}{d^{2-\beta}\varphi(D)} + \sum\limits_{d \le 2tD^{-2c}} \Lambda(d) \ll \frac{t}{\varphi(D)}. \notag
\end{align}

We are reduced to estimating the sum
\begin{equation}
J = \sum\limits_U \sum\limits_V \sum\limits_{V < d \le V+U^\beta} \Lambda(d) \left|Y(g,a\overline{d}, tV^{-1}) \right|. \notag
\end{equation}
Partitioning the variable $d$ according to the residue class $\modd{D}$ to which it belongs, a typical innersum over $d$ is
\begin{equation}
\sum\limits_{\substack{b = 1 \\ (b,D)=1}}^D \left| Y(g,a\overline{b}, tV^{-1} ) \right| \sum\limits_{\substack{V < d \le V+U^\beta \\ d \equiv b \modd{D}}} \Lambda(d) \notag
\end{equation}
\begin{equation}
\ll  \frac{U^\beta}{\varphi(D)} \sum\limits_{\substack{b=1 \\ (b,D)=1}}^D \left| Y(g,a\overline{b}, tV^{-1} ) \right|, \notag
\end{equation}
by a second application of the Brun-Titchmarsh theorem since, typically, $U^\beta \ge D^{\beta c}$ and we may choose $\beta$ so that $\beta c > 1$.

Bearing in mind that $a\overline{b}$ traverses a complete set of reduced residues $\modd{D}$ with $b$, an application of the Cauchy-Schwarz inequality shows this last bound to be
\begin{equation}
\ll \frac{U^\beta}{\varphi(D)} \left( \varphi(D) \sum\limits_{\substack{b=1 \\ (b,D)=1}}^D \left| Y(g,b,tV^{-1}) \right|^2 \right)^{1/2} \notag
\end{equation}
in turn, via the first waystation, Lemma \ref{JP_elliott_03_lem_I4},
\begin{equation}
\ll \frac{U^\beta}{\varphi(D)} \frac{tV^{-1}}{\log(tV^{-1})} \left( \frac{\log x}{\log D}\right)^\delta, \notag
\end{equation}
uniformly for $D^c \le U \le tD^{-2c}$, $U \le V \le 2U$, $D^c \le t \le y$, $x^\gamma \le y \le x$.

For the purposes of calculation it is convenient to replace this upper bound by
\begin{equation}
\ll \frac{1}{\varphi(D)} \left( \frac{\log x}{\log D} \right)^\delta \sum\limits_{V < d \le V+U^\beta} \frac{t}{d \log (td^{-1})} \notag
\end{equation}
to obtain the estimate
\begin{equation}
J \ll \frac{t}{\varphi(D)} \left( \frac{\log x}{\log D} \right)^\delta \sum\limits_{D^c < d \le 2tD^{-2c}} \frac{1}{d\log(td^{-1})}. \notag
\end{equation}
Since the function $(y\log(t/y))^{-1}$ is nonincreasing for $0 < y \le t/e$, the final sum is
\begin{equation}
\ll 1 + \int_{D^c}^{2tD^{-2c}} \frac{dw}{w\log(tw^{-1})} \ll \log\left( \frac{\log t}{\log D}\right). \notag
\end{equation}

Altogether,
\begin{equation}
Y(g\log, a, y) \ll \frac{y}{\varphi(D)} \left( \frac{\log x}{\log D} \right)^\delta \log\left( \frac{\log x}{\log D}\right), \notag
\end{equation}
uniformly for $x^\gamma \le y \le x$.

To remove the logarithm we note that
\begin{align}
& \sum\limits_{\substack{n \le y \\ n \equiv a \modd{D}}} g(n)\log(y/n) - \sum\limits_{j \in \mathcal{J}} \frac{\overline{\chi_j}(a)}{\varphi(D)} \sum\limits_{n \le y} g(n)\chi_j(n) \log(y/n) \notag \\
& \ll \sum\limits_{\substack{n \le y \\ n \equiv a \modd{D}}} \log(y/n) + \frac{1}{\varphi(D)} \sum\limits_{n \le y} \log(y/n) \ll \frac{y}{\varphi(D)}, \notag
\end{align}
uniformly for $1 \le D \le y/\log y$.

To complete the proof of Lemma \ref{JP_elliott_03_lem_I8} we set $\delta = \alpha/2$ and divide by $(\log y)^2$.

At this stage we remove the restrictions on the multiplicative function in the second waystation; first that it should vanish on the primes up to $D^c$.


\section{Truncated multiplicative functions.}


\makeatletter{}

\begin{lemma} \label{JP_elliott_03_lem_AA1}
Let $h$ be a real-valued multiplicative arithmetic function that for some $c_0$ satisfies $0 \le h(p^k) \le c_0^k$ on prime-powers and for each $\varepsilon > 0$, with an appropriate constant $c_1(\varepsilon)$, $h(n) \le c_1(\varepsilon)n^\varepsilon$ on the positive integers.  Let $0 < \beta < 1$.

Then
\begin{equation}
\sum\limits_{\substack{ n \le x \\ n \equiv a \modd{D}}} h(n) \ll \frac{x}{\varphi(D)\log x} \exp\left( \sum\limits_{\substack{ p \le x \\ (p,D)=1}} \frac{h(p)}{p} \right) \notag
\end{equation}
uniformly for $(a,D)=1$, $D \le x^\beta$, $x \ge 2$.
\end{lemma} 

\noindent \emph{Proof of Lemma \ref{JP_elliott_03_lem_AA1}}.  This is a particular case of a result of  Shiu \cite{shiu1980bruntitchmarsh}, that generalises the Brun-Titchmarsh theorem to nonnegative multiplicative functions.

\makeatletter{}

\begin{lemma} \label{JP_elliott_03_lem_AA2}
Let $2 \le w \le x$, $2 \le Y \le x$, and $0 < \beta < 1$.  Let $g$ be an exponentially multiplicative function with $0 \le g(p) \le 1$ if $p \le Y$, $g(p) = 0$ otherwise.

Then
\begin{equation}
\sum\limits_{\substack{ w < n \le x \\ n \equiv a \modd{D}}} g(n) \ll \frac{x}{\varphi(D)\log x} \exp\left( \sum\limits_{\substack{ p \le Y \\ (p,D)=1}} \frac{g(p)}{p} \right) \exp\left( - \frac{\log w}{\log Y} \right) \notag
\end{equation}
uniformly for $(a,D)=1$, $D \le x^\beta$.
\end{lemma} 

\makeatletter{}

\noindent \emph{Proof of Lemma \ref{JP_elliott_03_lem_AA2}}.  Set $\theta = 1/\log Y$.  The sum, $S$, to be estimated, does not exceed
\begin{equation}
w^{-\theta} \sum\limits_{\substack{ n \le x \\ n \equiv a \modd{D}}} g(n) n^\theta. \notag
\end{equation}
Set $h(n) = g(n)n^\theta$.  Here $h(p^k) = \left( g(p)\exp(\theta \log p) \right)^k / k! \le e^k/k!$ if $p \le Y$, and is otherwise zero.  Thus $h(p^k) \le \exp(e)$ for all prime-powers $p^k$.  In particular, $h(n) \le \exp\left( e\omega(n) \right)$, where $\omega(n)$ counts the number of distinct prime-divisors of $n$, hence $\le \exp( c_2 \log n/\log\log n) < n^\varepsilon$ for all $n$ sufficiently large.

After an application of Lemma \ref{JP_elliott_03_lem_AA1},
\begin{equation}
S \ll \frac{w^{-\theta}x}{\varphi(D)\log x} \exp\left( \sum\limits_{\substack{ p \le Y \\ (p,D)=1}} g(p)p^{\theta-1} \right). \notag
\end{equation}
Appealing to the bound $e^t - 1 \le te^t$, valid for all nonnegative $t$, the sum in the exponential does not exceed
\begin{equation}
\sum\limits_{\substack{ p \le Y \\ (p,D)=1}} g(p) p^{-1} + e\theta \sum\limits_{p \le Y} g(p) p^{-1} \log p, \notag
\end{equation}
with a second term that is $\ll \theta \log Y \ll 1$.

Since $w^{-\theta} = \exp( -\log w/\log Y)$, the lemma is established. 

\makeatletter{}

\begin{lemma} \label{JP_elliott_03_lem_AA3}
Let $2 \le w \le x$, $2 \le Y \le x$.  Let $g$ be a multiplicative function with values in the real interval, $0 \le g(n) \le 1$, vanishing on the primes $p > Y$.

Then
\begin{equation}
\sum\limits_{\substack{ w < n \le x \\ n \equiv a \modd{D}}} g(n) \ll \frac{x}{\varphi(D)\log w} \prod\limits_{\substack{ p \le Y \\ (p,D)=1}} \left( 1 + \frac{g(p)}{p} \right) \exp\left( - \frac{\log w}{10\log Y} \right) \notag
\end{equation}
uniformly for $(a,D)=1$, $D^{11} \le w$.
\end{lemma} 

 \emph{Remark}.  Although susceptible of improvement, the uniformities in Lemma \ref{JP_elliott_03_lem_AA3} are adequate to our present requirements.

\makeatletter{}

\noindent \emph{Proof of Lemma \ref{JP_elliott_03_lem_AA3}}.  Express $g$ as a convolution of multiplicative functions $\ell * r$, where $\ell(p^k) = g(p)^k / k!$, $k \ge 1$.  An examination of Euler products shows that $r(p) = 0$ and
\begin{equation}
r(p^k) = \sum\limits_{j=0}^k \frac{1}{j!} \left( -g(p) \right)^j g(p^{k-j}), \quad k \ge 2. \notag
\end{equation}
Hence $|r(p^2)| \le 1$ and $|r(p^k)| \le e$ generally.  It is convenient to note that
\begin{equation}
\sum\limits_v |r(v)| v^{-1}  \le \prod\limits_p \left( 1 + ep^{-2}(1-p^{-1})^{-1} \right) \notag
\end{equation}
and absolutely bounded.

Let $0 < \delta < 1$.  We represent
\begin{equation}
\sum\limits_{\substack{ w < n \le x \\ n \equiv a \modd{D}}} g(n) = \sum\limits_{v \le xw^{-\delta}} r(v) \sum\limits_{\substack{ w^\delta < u \le xv^{-1} \\ u \equiv a\overline{v} \modd{D}}} \ell(u) + \sum\limits_{u \le w^\delta} \ell(u) \sum\limits_{\substack{ wu^{-1} < v \le xu^{-1} \\ v \equiv a \overline{u} \modd{D}}} r(v), \notag
\end{equation}
$ = \Sigma_1 + \Sigma_2$, say.

Applied to its innersum, Lemma \ref{JP_elliott_03_lem_AA2} assures that
\begin{equation}
\Sigma_1 \ll \frac{x}{\varphi(D)\log w} \exp\left( \sum\limits_{\substack{ p \le Y \\ (p,D)=1}} \frac{g(p)}{p} - \frac{\delta \log w}{\log Y} \right) \notag
\end{equation}
uniformly for $D \le (xv^{-1})^\beta$, hence for $D \le w^{\delta\beta}$.

If a prime $p$ exactly divides an integer $v$ in $\Sigma_2$, then $r(v) = 0$.  Thus $v = m^2 t$ where $p \mid t$ implies $p \mid m$.  In particular, $t \le m$, so $m \ge (wu^{-1})^{1/3}$.  Moreover, $m \le x^{1/2}$.  The innersum in $\Sigma_2$ is
\begin{align}
& \ll \sum\limits_{(wu^{-1})^{1/3} < m \le x^{1/2}} e^{\omega(m)} \sum\limits_{\substack{ t \le xm^{-2} \\ t \equiv a\overline{m}^2 \modd{D}}} 1 \notag \\
& \ll \sum\limits_{(wu^{-1})^{1/3} < m \le x^{1/2}} e^{\omega(m)} \left( \frac{x}{m^2 D} + 1 \right). \notag
\end{align}
Estimating simply,
\begin{equation}
\sum\limits_{m \le y} e^{\omega(m)} \le y \sum\limits_{m \le y} e^{\omega(m)} m^{-1} \le y \prod\limits_{p \le y} \left( 1+ep^{-1} + O(p^{-2}) \right) \ll y(\log y)^e, \notag
\end{equation}
and with an integration by parts,
\begin{equation}
\sum\limits_{m > y} e^{\omega(m)} m^{-2} \ll y^{-1} (\log y)^e, \quad y \ge 2. \notag
\end{equation}
Hence
\begin{align}
\Sigma_2 & \ll \sum\limits_{u \le w^\delta} \ell(u) \left( \frac{x}{D} \left(\frac{u}{w} \right)^{1/3} (\log w)^e + x^{1/2}(\log x)^e \right) \notag \\
& \ll \sum\limits_{u \le w^\delta} \ell(u)u^{-1} \left( xD^{-1} w^{4\delta/3 - 1/3}(\log w)^e + w^\delta x^{1/2}(\log x)^e \right). \notag
\end{align}

We choose $\delta = 1/10$.  The series involving $\ell$ is
\begin{equation}
\le \prod\limits_{\substack{ p \le Y \\ (p,D)=1}} \left( 1+ \sum\limits_{k=1}^\infty g(p)^k/k! \right) = \exp\left( \sum\limits_{\substack{ p \le Y \\ (p,D)=1}} \frac{g(p)}{p} \right), \notag
\end{equation}
the coefficient of $xD^{-1}$ is
\begin{equation}
w^{-1/5} (\log w)^e \ll \exp(-\log w/10\log Y)(\log w)^{-1} \notag
\end{equation}
and the restriction on the size of $D$ guarantees the term involving $x^{1/2}$ to be negligible.

This completes the proof of Lemma \ref{JP_elliott_03_lem_AA3}.

It is useful to note the following

\noindent {\bf Corollary to Lemma \ref{JP_elliott_03_lem_AA3}.}  \emph{Let $2 \le w \le x$, $2 \le Y \le x$.  Let $g$ be a multiplicative function with values in the complex unit disc, vanishing on the primes $p > Y$.}

\emph{Then
\begin{equation}
Y(g,a,x) - Y(g,a,y) \ll \frac{x}{\varphi(D)\log w} \prod\limits_{\substack{p \le Y \\ (p,D)=1}} \left(1+ \frac{|g(p)|}{p} \right) \exp\left( -\frac{\log w}{10\log Y} \right) \notag
\end{equation}
uniformly for $w \le y \le x$, $D^{11} \le w$.}

\noindent \emph{Proof of Corollary to Lemma \ref{JP_elliott_03_lem_AA3}}.  The difference to be estimated is
\begin{equation}
\ll \sum\limits_{\substack{ w < n \le x \\ n \equiv a \modd{D}}} |g(n)| + \frac{1}{\varphi(D)} \sum\limits_{\substack{ w < n \le x \\ (n,D)=1}} |g(n)|. \notag
\end{equation}
We may apply Lemma \ref{JP_elliott_03_lem_AA3} to the first sum, directly, and to the second sum after partitioning the range of the variable $n$ into reduced residue classes $\modd{D}$.


\section{Proof of Theorem \ref{JP_elliott_03_thm_01}: Completion}


Express $g$ as a convolution of multiplicative functions $f*h$, where $f(p^k) = g(p^k)$ if $p > D^c$ and is zero otherwise, $h$ coincides with $g$ on the powers of $p \le D^c$ and is zero otherwise.  There is a corresponding decomposition
\begin{equation}
Y(g,a,y) = \sum\limits_{v \le y^{1/2}} h(v)Y(f,a\overline{v}, yv^{-1}) + \sum\limits_{u \le y^{1/2}} f(u) \left( Y(h, a\overline{u}, yu^{-1}) - Y(h,a\overline{u}, y^{1/2}) \right). \notag
\end{equation}
From the second Waystation the sum over the variable $v$ is
\begin{equation}
\ll \sum\limits_{v \le y^{1/2}} |h(v)| \frac{y}{v\varphi(D)\log y} \left( \frac{\log y}{\log D}\right)^\alpha \notag
\end{equation}
where
\begin{equation}
\sum\limits_{v \le y^{1/2}} |h(v)| v^{-1} \ll \prod\limits_{p \le D^c} \left( 1+ |g(p)| p^{-1} \right) \ll \prod\limits_{p \le D} \left( 1 + |g(p)|p^{-1} \right). \notag
\end{equation}
In view of the Corollary to Lemma \ref{JP_elliott_03_lem_AA3}, the sum over the variable $u$ is
\begin{equation}
\ll \sum\limits_{u \le y^{1/2}} \frac{|f(u)|}{u} \frac{y}{\varphi(D)\log y} \prod\limits_{p \le D} \left( 1+ \frac{|g(p)|}{p} \right) \exp\left( -\frac{\log y}{20c\log D} \right). \notag
\end{equation}
Here
\begin{equation}
\sum\limits_{u \le  y^{1/2}} \frac{|f(u)|}{u} \ll \prod\limits_{D^c < p \le y^{1/2}} \left(1 + \frac{|g(p)|}{p} \right) \ll \frac{\log y}{\log D}, \notag
\end{equation}
an amount which is absorbed by the previous exponential factor.

It remains to remove the restriction that $g$ be completely multiplicative.  Since the argument runs along familiar lines, we indicate only salient details.

Given a convolution representation $g = \ell * r$ with $\ell(p^k) = g(p)^k$, $k \ge 1$, examination of Euler products shows that $r(p^k) = g(p^k) - g(p)g(p^{k-1})$, $k \ge 1$, in particular that $r(p)=0$.

We employ the above decomposition of $Y(g,a,y)$, the roles of $f$, $h$ played by $\ell$, $r$ respectively.  To the corresponding sum over $v \le y^{1/2}$ we apply the version of Theorem \ref{JP_elliott_03_thm_01} obtained so far.  The sum over $u \le y^{1/2}$ we treat in the manner of Lemma \ref{JP_elliott_03_lem_AA3}, the role of the function $e^{\omega(m)}$ played by $2^{\omega(m)}$.

This completes the proof of Theorem \ref{JP_elliott_03_thm_01}.

Note that, in accordance with the second of the remarks following the definition of exceptional characters, Theorem \ref{JP_elliott_03_thm_01} applies simultaneously to the functions $g$ and $g\mu$.


\section{Taxonomy of exceptional characters} \label{JP_elliott_03_sec_taxonomy}


In this section we address the exceptional characters that appear in Theorem \ref{JP_elliott_03_thm_01}.  An overview is that without hypothesis at most one Dirichlet character $\modd{D}$ can be near to a given multiplicative function $g$ with values in the complex unit disc.  Unless $g$ has slender support, all exceptional characters in Theorem \ref{JP_elliott_03_thm_01} are obtained by braiding a character close to $g$ with characters of order bounded independently of the modulus $D$.

To motivate the methodology we apply it to derive the following result.

\makeatletter{}

\begin{theorem} \label{JP_elliott_03_thm_A}
Let $c$, $c_1$, $\gamma$, $x$ be positive real numbers, $D$ an integer, $0 < c < 1$, $0 < \gamma < 1$, $1 \le D \le x$.

Let $g$ be a multiplicative function, with values in the complex unit disc, that satisfies
\begin{equation}
\sum\limits_{w < p \le x} \left(|g(p)| - c\right) p^{-1} \ge -c_1, \quad D \le w \le x, \notag
\end{equation}
on the primes.

Then there are nonprincipal Dirichlet characters $\chi_j \modd{D}$, their number bounded in terms of $c$ alone, such that uniformly for $(a,D)=1$, $x^\gamma \le y \le x$,
\begin{align}
\sum\limits_{\substack{n \le y \\ n \equiv a \modd{D}}} g(n) = \frac{1}{\varphi(D)} & \sum\limits_{\substack{n \le y \\ (n,D)=1}} g(n) + \sum\limits_j \frac{\overline{\chi_j}(a)}{\varphi(D)} \sum\limits_{n \le y} g(n)\chi_j(n) \notag \\
& + O\Biggl( \frac{y}{\varphi(D)\log y} \prod\limits_{\substack{p \le y \\ (p,D)=1}} \left(1+ \frac{|g(p)|}{p} \right) \left(\frac{\log D}{\log y}\right)^\eta \Biggr) \notag
\end{align}
with $\eta = 10^{-3}c^4(c+1)^{-1}$ and the order of every product $\chi_j \overline{\chi_r}$, $j \ne r$, of exceptional characters bounded by $10c^{-1}$.

If for some positive integer $k$
\begin{equation}
\sum\limits_{\substack{D < p \le x \\ g(p)^k \in \mathbb{R}}} \left(|g(p)| - c\right) p^{-1} \ge -c_1 \notag
\end{equation}
and the exponent $\eta$ is replaced by $k^{-2}\eta$, then the order of each exceptional character is at most $20kc^{-1}$.

Moreover, with $\left( 10^{-1} c k^{-1} \right)^2 \eta$ in place of $\eta$ and $Z = \exp\left( \log D(\log x/\log D)^{1/30} \right)$, each exceptional sum satisfies the uniform bound
\begin{equation}
\sum\limits_{n \le y} g(n) \chi_j(n) \ll \frac{y}{\log y} \prod\limits_{\substack{p \le y \\ (p,D)=1}} \left(1+\frac{|g(p)|}{p} \right) \exp\left(-\frac{c}{c+1} \sum\limits_{D < p \le Z} \frac{|g(p)| - \Re g(p)\chi_j(p)}{p} \right), \notag
\end{equation}
$x^\gamma \le y \le x$, the implied constant depending at most upon $c$, $c_1$, $\gamma$ and $k$.
\end{theorem}


\emph{Remarks}.  Integrating by parts, the lower bound hypothesis on $|g(p)|$ follows directly from a uniform lower bound
\begin{equation}
\sum\limits_{p \le y} |g(p)| p^{-1} \log p \ge c\log y - c_2, \quad y \ge 2. \notag
\end{equation}

Note that
\begin{equation}
\sum\limits_{D < p \le Z} p^{-1} = \frac{1}{30} \sum\limits_{D < p \le x} p^{-1} + O(1), \quad 1 \le D \le x. \notag
\end{equation}

Theorem \ref{JP_elliott_03_thm_A} aims to maintain uniformity in the modulus $D$ whilst encompassing multiplicative functions $g$ with limited and scattered support and reducing the orders of the exceptional characters to a range within the purview of standard reciprocity laws.

The constants appearing in Theorem \ref{JP_elliott_03_thm_A} may be considerably varied.  In practice the lower bound on $|g(p)|$ with a possibly small value of $c$ serves to enable the application of Theorem \ref{JP_elliott_04_thm_02} to Theorem \ref{JP_elliott_03_thm_01} and support a taxonomy of the exceptional characters according to the size of the sum $\sum_{D < p \le x} |g(p)|p^{-1}$ alone, as is shown in Theorem \ref{JP_elliott_03_thm_B}.

For $T \ge D \ge 2$, denote by $\Delta(T)$ the function $\log(\log T/\log D)+c$ that appears as the upper bound in Lemma \ref{JP_elliott_03_lem_I3}.


\makeatletter{}

\begin{lemma}\label{JP_elliott_03_lem_A1}
Let $\chi$ be a Dirichlet character $\modd{D}$, $D \ge 2$, $x$, $\delta$, $t$ real, $x \ge D$, $0 < \delta \le 1$, $|t| \le T$, $h$ a real-valued function on the primes in the interval $(D,x]$ for which $0 \le h(p) \le 1$ and
\begin{equation}
\sum\limits_{D < p \le x} h(p) p^{-1} \left| 1-\chi(p)p^{it} \right|^2 \le \delta L, \notag
\end{equation}
where $L = \sum_{D < p \le x} p^{-1}$.  

Then
\begin{center}
either $\displaystyle \sum\limits_{D < p \le x} h(p)p^{-1} \le 4\delta^{1/3}L + \Delta\left(\delta^{-1/3}T\right)$ or the order of $\chi$ is less than $2\delta^{-1/3}$.
\end{center}
\end{lemma} 


\emph{Remark}.  In application to Theorem \ref{JP_elliott_03_thm_A}, $\delta$ will be fixed at a positive value independent of $D$, $x$;  $T$ will be chosen of the form $\max\left(D, (\log x/\log D)^d\right)$ for some positive constant $d$, also independent of $D$, $x$.  As a consequence, the term $\Delta \left( \delta^{-1/3}T\right)$ in the above conclusion will not exceed $\log L + O(1)$, effectively negligible in comparison with $L$.


\makeatletter{}

\noindent \emph{Proof of Lemma \ref{JP_elliott_03_lem_A1}}.  We begin with the Fej\'er kernel
\begin{equation}
\frac{1}{N}\left(\frac{\sin \pi N\theta}{\sin \pi\theta}\right)^2 = \sum\limits_{m = -(N-1)}^{N-1} \left(1-\frac{|m|}{N}\right)e^{2\pi i m \theta}, \quad \theta \in \mathbb{R}. \notag
\end{equation}
It is useful to note that if $\| \theta \|$ denotes the distance of $\theta$ to a nearest integer, then
\begin{equation}
2\|\theta\| \le \left|\sin \pi \theta\right| \le \pi \|\theta\|. \notag
\end{equation}
If, moreover, $2N\|\theta\| \le 1$, then $\|N\theta\| = N\|\theta\|$ and the Fej\'er kernel is at least $4N\pi^{-2}$.  With $2\pi\gamma_p = \arg\chi(p) + t\log p$ we see that
\begin{align}
4N\pi^{-2} \sum\limits_{\substack{D < p \le x \\ 2N\|\gamma_p\| \le 1}} \frac{h(p)}{p} & \le \sum\limits_{D < p \le x}\frac{1}{N}\left(\frac{\sin \pi N\gamma_p}{\sin\pi\gamma_p}\right)^2 \frac{1}{p} \notag \\
& = \sum\limits_{m=-N+1}^{N-1} \left(1-\frac{|m|}{N}\right)\text{Re} \sum\limits_{D < p \le x} \left(\chi(p)p^{it}\right)^m \frac{1}{p}. \notag
\end{align}
If the order of $\chi$ is at least $N$, then each of the innermost character sums with $m \ne 0$ is by Lemma \ref{JP_elliott_03_lem_I3} at most $\Delta(NT)$.  

However, $\left|1-\chi(p)p^{it}\right| = \left|1-e^{2\pi i \gamma_p}\right| = 2\left|\sin\pi\gamma_p\right| \ge 4\|\gamma_p\|$ so that
\begin{equation}
\sum\limits_{\substack{D<p \le x \\ 2N\|\gamma_p\| > 1}} \frac{h(p)}{p} < \frac{N^2}{4} \sum\limits_{D < p \le x} \frac{h(p)}{p} \left| 1 - \chi(p)p^{it} \right|^2 \le \frac{N^2\delta L}{4}. \notag
\end{equation}
Altogether,
\begin{equation}
\sum\limits_{D<p \le x} \frac{h(p)}{p} \le \left( \frac{N^2\delta}{4} + \frac{\pi^2}{4N} \right)L + \Delta(NT). \notag
\end{equation}
Choosing, for simplicity of exposition, $N=[2\delta^{-1/3}]$, so that $N>2\delta^{-1/3}-1 \ge \delta^{-1/3}$, we obtain the upper bound $4\delta^{1/3}L +\Delta\left(\delta^{-1/3}N\right)$ and complete the proof. 


Employing the lower bound $\sin \theta \ge \theta - \theta^3/6$, valid for $0 \le \theta \le \frac{\pi}{2}$, and the consequent lower bound
\begin{equation}
\frac{1}{r}\left(\frac{\sin \pi r\theta}{\sin \pi \theta}\right)^2 \ge r\left(1-\frac{(\pi r \|\theta\|)^2}{6}\right) \notag
\end{equation}
for the Fej\'er kernel, valid if $2r\|\theta\| \le 1$, a similar argument with $r$ in place of $N$ and a division into cases according to whether $2r\| \gamma_p \| \le \left(3\pi^{-2} r^3\delta \right)^{1/4}$ delivers the following variant of Lemma \ref{JP_elliott_03_lem_A1}.


\makeatletter{}

\begin{lemma}\label{JP_elliott_03_lem_A2}
Under the hypotheses of Lemma \ref{JP_elliott_03_lem_A1}, if $r^3\delta \le 1$ and $\chi$ has order at least $r$, $r \ge 2$, then
\begin{equation}
\sum\limits_{D<p \le x} \frac{h(p)}{p} \le \left(1 + (r^3\delta)^{1/2}\right)\frac{L}{r} + \Delta\left((r-1)T\right). \notag
\end{equation}
\end{lemma} 


If $D$ is a prime, $D \equiv 1 \modd{r}$ and $S$ is the set of primes in $(D,x]$ that are $r^{th}$-powers$\modd{D}$, then for every character of order $r$: $\sum_{p \in S} p^{-1} \left|1-\chi(p)\right|^2 = 0$ whilst, for a fixed $D$, Dirichlet's Theorem on primes in arithmetic progressions shows that $\sum_{p \in S}p^{-1} = r^{-1}L + O(1)$ as $x \to \infty$.  The upper bound in Lemma \ref{JP_elliott_03_lem_A2} cannot be appreciably improved.

Under the hypotheses of Lemma \ref{JP_elliott_03_lem_A1} there is control on the size of $t$, too.  To this end we employ an analogue of Lemma \ref{JP_elliott_03_lem_I3}.


\makeatletter{}

\begin{lemma} \label{JP_elliott_03_lem_A3}
With a certain real $c_0$,
\begin{equation}
\Re \sum\limits_{y < p \le x} p^{-1-it} \le \begin{cases} 2\log\log\left(2+|t|\right) + c_0, & |t|>(\log y)^{-1} \\
-\log\left(|t|\log y\right) + c_0, & (\log x)^{-1} < |t| \le (\log y)^{-1}, \\ \log(\log x/\log y) + c_0, & |t| \le (\log x)^{-1}, \end{cases} \notag
\end{equation}
uniformly for $x \ge y \ge 2$, $t$ real.
\end{lemma} 

\makeatletter{}

\noindent \emph{Proof of Lemma \ref{JP_elliott_03_lem_A3}}.  Applying the estimates of Chebyshev as in Lemma \ref{JP_elliott_03_lem_I3},
\begin{equation}
\exp\left( \Re \sum\limits_{y < p \le x} p^{-1-it} \right) \notag
\end{equation}
lies between positive absolute constant multiples of
\begin{equation}
\left| \zeta \left(1+(\log x)^{-1} + it \right) \left( \zeta\left(1+(\log y)^{-1} + it \right) \right)^{-1} \right|. \notag
\end{equation}

The first of the three inequalities in the lemma follows from the classical bounds $\left(\log(2+|t|)\right)^{-1} \ll |\zeta(s)| \ll \log(2+|t|)$, valid in the notched half-plane $\sigma \ge 1$, $|s-1| \ge c_1 > 0$, c.f. \cite{titchmarsh1986zeta}, Theorem 3.5.  The second and third inequalities follow from the Laurent expansion $\zeta(s) = (s-1)^{-1} + \cdots$ around the simple pole of $\zeta(s)$ at $s=1$. 


Denote the upper bound in Lemma \ref{JP_elliott_03_lem_A3} by $\psi(t)$ and set $L_1 = \sum_{y < p \le x} p^{-1}$.


\makeatletter{}

\begin{lemma} \label{JP_elliott_03_lem_A4}
If for some $\delta$, $0 < \delta < 1$, $t$ real with $h(p)$ in the interval $[0,1]$
\begin{equation}
\sum\limits_{y < p \le x} h(p)p^{-1} \left|1-p^{it}\right|^2  \le \delta L_1, \notag
\end{equation}
then, with the same uniformities
\begin{equation}
\sum\limits_{y < p \le x} h(p)p^{-1} \le 4\delta^{1/3}L_1 + 3\psi\left(2\delta^{-1/3}t\right) + 3\psi(t). \notag
\end{equation}
\end{lemma} 


\noindent \emph{Proof of Lemma \ref{JP_elliott_03_lem_A4}}.  The argument follows that for Lemma \ref{JP_elliott_03_lem_A1} with appeal to Lemma \ref{JP_elliott_03_lem_A3} in place of that to Lemma \ref{JP_elliott_03_lem_I3}.


Before proceeding to the proof of Theorem \ref{JP_elliott_03_thm_A}, let $D$ be an integer, $x$, $T$ real, $2 \le D \le x$, $D \le T$ and again $\Delta(T)$ the function $\log(\log T/\log D)+c$ that appears as the upper bound in Lemma \ref{JP_elliott_03_lem_I3}.  In the following proof $T$ is chosen to be $\max\left(D, (\log x/\log D)^4\right)$, so that the remark made following the statement of Lemma \ref{JP_elliott_03_lem_A1} comes into force.

Assuming that $L = \sum_{D < p \le x} p^{-1} > 0$, on the equivalence classes of multiplicative functions, with values in the complex unit disc, and that coincide on the primes in the interval $(D,x]$, define the metric
\begin{equation}
\rho(g_1, g_2) = \left( \frac{1}{4L} \sum\limits_{D < p \le x} \frac{1}{p} \left| g_1(p) - g_2(p) \right|^2 \right)^{1/2} \ge 0. \notag
\end{equation}
Effectively, multiplicative functions with values in the complex unit disc belong to a ball of radius 1 whose centre we may choose to be the arithmetic function that is identically 1.

For distinct Dirichlet characters $\chi_j \modd{D}$, real $t_j$, $|t_j| \le T$, define generalised characters $\chi_{j, t_j}$ by $p \to \chi_j(p)p^{it_j}$, $j = 1, 2$.  An application of Lemma \ref{JP_elliott_03_lem_I3} shows that
\begin{equation}
\rho(\chi_{1,t_1}, \chi_{2,t_2} ) = \left( \frac{1}{2L} \sum\limits_{D < p \le x} \frac{1}{p} \left(1 - \text{Re}\, \chi_{1,t_1}\overline{\chi_{2,t_2}}(p) \right) \right)^{1/2} \notag
\end{equation}
\begin{equation}
\ge \left(\frac{1}{2} - \frac{\Delta(2T)}{2L} \right)^{1/2} = \frac{1}{\sqrt{2}} + O\left(\frac{\log L}{L}\right) \notag
\end{equation}
and, after the triangle inequality,
\begin{equation}
\max\limits_{j=1,2} \rho\left(\chi_{j,t_j}, g\right) \ge \frac{1}{2\sqrt{2}} + O\left( \frac{\log L}{L} \right). \notag
\end{equation}

To this extent, at most one generalised character can be near to a given multiplicative function, $g$.

\makeatletter{}

\noindent \emph{Proof of Theorem \ref{JP_elliott_03_thm_A}}.  For primes on which $g$ does not vanish, let $g(p) = |g(p)|e^{i\theta_p}$; otherwise set $\theta_p=0$.

We apply Theorem \ref{JP_elliott_03_thm_01} with $\alpha$ a positive value not exceeding $c-\eta$.  Under the lower bound hypothesis on $|g(p)|$ this ensures that the error term in Theorem \ref{JP_elliott_03_thm_01} falls within that of Theorem \ref{JP_elliott_03_thm_A}.

If, for an exceptional character $\chi$, real $\delta$, $0 < \delta < 1$,
\begin{equation}
\min\limits_{|t| \le T/2} \sum\limits_{D < p \le x} |g(p)| \left| 1-e^{i\theta_p} \chi(p) p^{it} \right|^2 p^{-1} > \delta L/4, \notag
\end{equation}
then with $Y=D$ the function $\lambda$ in Theorem \ref{JP_elliott_04_thm_02} exceeds $\delta L/8$.  With $\delta = (c/5)^3$ the corresponding sum over the $g\chi(n)$ also falls within the error term of Theorem \ref{JP_elliott_03_thm_A}.

For any pair $\chi_j$, $\chi_r$ of the remaining exceptional characters an application of the Cauchy-Schwarz inequality shows that for certain real $t_j$, $t_r$ in the interval $[-T/2, T/2]$
\begin{equation}
\sum\limits_{D < p \le x} |g(p)| \left|1-\chi_j \overline{\chi_r}(p)p^{i(t_j-t_r)} \right|^2 p^{-1} \le \delta L. \notag
\end{equation}
Since the inequality $cL - c_1 \le 4\delta^{1/3}L + \Delta\left( \delta^{-1/3}T\right)$ with $4\delta^{1/3} = 4c/5$ fails provided $D$ does not exceed a sufficiently small fixed power of $x$, which we may assume, an application of Lemma \ref{JP_elliott_03_lem_A1} guarantees the order of $\chi_j\overline{\chi_r}$ not to exceed $2\delta^{-1/3}$, i.e., $10c^{-1}$.

Suppose, further, that
\begin{equation}
\sum\limits_{\substack{D<p \le x \\ g(p)^k \in \mathbb{R}}} \left(|g(p)| - c\right)p^{-1} \ge -c_1 \notag
\end{equation}
and that for some real $t$, $|t| \le T/(2k)$,
\begin{equation}
\sum\limits_{\substack{D < p \le x \\ g(p)^k \in \mathbb{R}}} |g(p)| \left|1-e^{i\theta_p} \chi(p)p^{it} \right|^2 p^{-1} \le \left(4k^2\right)^{-1}\delta L. \notag
\end{equation}
From the inequality $|1-z^k| \le k|1-z|$, valid for $z$ in the complex unit disc,
\begin{equation}
\sum\limits_{\substack{D < p \le x \\ g(p)^k \in \mathbb{R}}} |g(p)| \left|1-\left(\chi(p)p^{it} \right)^{2k} \right|^2 p^{-1} \le \delta L \notag
\end{equation}
and an application of Lemma \ref{JP_elliott_03_lem_A1} shows $\chi^{2k}$ to have an order not exceeding $10c^{-1}$.

In particular, provided we replace $\eta$ by $k^{-2}\eta$, we may assume that every exceptional character has order not exceeding $20kc^{-1}$.


Let $m$ be the order of a typical character $\chi^{2k}$.  At the expense of replacing $k^{-2}\eta$ by $(mk)^{-2}\eta$, i.e., of replacing $\delta$ by $m^{-2}\delta$, we may assume that for some $t$ in the shorter range
\begin{equation}
\sum\limits_{\substack{D < p \le x \\ g(p)^k \in \mathbb{R}}} |g(p)| \left|1-p^{2kmit} \right|^2 p^{-1} \le \delta L \notag
\end{equation}
and may apply Lemma \ref{JP_elliott_03_lem_A4} with $y=D$ to conclude that
\begin{equation}
cL-c_1 \le 4\delta^{1/3}L + 3\psi\left(2\delta^{-1/3}t \right) + 3\psi(t). \notag
\end{equation}
If $|t|\log D > 1$ then the second and third terms in this bound do not exceed $6\log L + O(1)$ which, for $\log x/\log D$ sufficiently large, is untenable.  If $(\log x)^{-1} < |t| \le (\log D)^{-1}$ then $cL/5 \le -6\log\left(|t|\log D\right) + O(1)$, hence $|t| \log Z \ll 1$.

In particular, replacing $p^{it}$ by 1 in the corresponding sum
\begin{equation}
\sum\limits_{p \le Z} \left(|g(p)| - \Re g(p)\chi(p) p^{it} \right)p^{-1} \notag
\end{equation}
introduces an error of
\begin{equation}
\ll \sum\limits_{p \le Z} \left|1-p^{it}\right| p^{-1} \ll |t| \sum\limits_{p \le Z} p^{-1} \log p \ll |t| \log Z \ll 1. \notag
\end{equation}

We may argue similarly for every $t$ such that
\begin{equation}
\sum\limits_{D < p \le x} \left(|g(p)| - \Re g(p)\chi(p)p^{it} \right) p^{-1} \le \left(20kc^{-1}\right)^{-2} \delta L \notag
\end{equation}
and the final assertion of Theorem \ref{JP_elliott_03_thm_A} rapidly follows. 


For convenience of application we state a version of Theorem \ref{JP_elliott_03_thm_A} that is obtained by employing the refined Lemma \ref{JP_elliott_03_lem_A2} in place of Lemma \ref{JP_elliott_03_lem_A1}.


\makeatletter{}

\begin{theorem} \label{JP_elliott_03_thm_B}
Let $c$, $c_1$, $\gamma$, $\varepsilon$, $x$ be positive real numbers, $D$ an integer, $0 < c < 1$, $0 < \gamma < 1$, $0 < \varepsilon < 1$, $1 \le D \le x$.

Let $g$ be a multiplicative function, with values in the complex unit disc, that satisfies
\begin{equation}
\sum\limits_{w < p \le x} \left(|g(p)| - c\right) p^{-1} \ge -c_1, \quad D \le w \le x, \notag
\end{equation}
on the primes and, for a positive integer $r$,
\begin{equation}
\sum\limits_{D < p \le x} |g(p)| p^{-1} > \left(\frac{1}{r} + \varepsilon \right) \sum\limits_{D < p \le x} p^{-1} + c_0, \notag
\end{equation}
where $c_0$ is the constant appearing in Lemma \ref{JP_elliott_03_lem_A2}.

Then there are nonprincipal Dirichlet characters $\chi_j \modd{D}$, their number bounded in terms of $c$ alone, such that uniformly for $(a,D)=1$, $x^\gamma \le y \le x$,
\begin{align}
\sum\limits_{\substack{n \le y \\ n \equiv a \modd{D}}} g(n) = \frac{1}{\varphi(D)} & \sum\limits_{\substack{n \le y \\ (n,D)=1}} g(n) + \sum\limits_j \frac{\overline{\chi_j}(a)}{\varphi(D)} \sum\limits_{n \le y} g(n)\chi_j(n) \notag \\
& + O\Biggl(\frac{y}{\varphi(D)\log y} \prod\limits_{\substack{p \le y \\ (p,D)=1}} \left(1 + \frac{|g(p)|}{p} \right) \left( \frac{\log D}{\log y} \right)^\eta \Biggr) \notag
\end{align}
with $\eta = c\varepsilon^2\left(2r^3(c+1)\right)^{-1}$ and the order of every product $\chi_j \overline{\chi_k}$, $j \ne k$, of exceptional characters less than $r$.

Moreover, if $g$ is real and we replace $\eta$ by $\eta/2$, then every exceptional character has a square of order less than $r$.
\end{theorem} 


\begin{example} \label{JP_elliott_03_exa_01}
If $g = \mu$, the M\"obius function, then with $r=2$, $2\sqrt{2} \varepsilon=1$, $c=1$, $\eta = 2^{-9}$, there can be at most one exceptional character, and that real.
\end{example}


\begin{example} \label{JP_elliott_03_exa_02}
As demonstrated by Landau over a century ago, information on the value distribution of M\"obius' function informs the distribution of prime numbers.  As a second example we give a proof of Linnik's theorem that for a positive constant $c$, the least prime in each reduced residue class $\modd{D}$ does not exceed $D^c$ in size.

In the interest of brevity we employ only modest values for various parameters.

\makeatletter{}

\begin{lemma} \label{JP_elliott_03_lem_BB13}
Let $0 < \gamma < 1$, $2 \le D \le N^\gamma$.  Then there exists a nonprincipal real character $\chi \modd{D}$ such that with $\tau = 2^{-10}$
\begin{equation}
\sum\limits_{\substack{ D < n \le x \\ n \equiv a \modd{D}}} \Lambda(n) = \frac{1}{\varphi(D)} \sum\limits_{\substack{D < n \le x \\ (n,D)=1}} \Lambda(n) + \frac{\chi(a)}{\varphi(D)} \sum\limits_{D < n \le x} \Lambda(n)\chi(n) + O\left( \frac{x}{\varphi(D)} \left( \frac{\log D}{\log x} \right)^\tau \right) \notag
\end{equation}
uniformly for $N^\gamma \le x \le N$, $(a,D)=1$.
\end{lemma} 

An integration by parts then delivers the

\begin{corollary}
Under the hypotheses of Lemma \ref{JP_elliott_03_lem_BB13},
\begin{equation}
\varphi(D) \sum\limits_{\substack{ N^\gamma < p \le N \\ p \equiv a \modd{D}}} p^{-1} = \sum\limits_{N^\gamma < p \le N} p^{-1} + \chi(a) \sum\limits_{N^\gamma < p \le N} p^{-1} \chi(p) + O\left((\log D/\log N)^\tau \right), \notag
\end{equation}
the implied constant depending at most upon $\gamma$.
\end{corollary}

\makeatletter{}

\noindent \emph{Proof of Lemma \ref{JP_elliott_03_lem_BB13}}.  Define the completely multiplicative function $g$ by $g(p) = 1$ if $p>D$, $g(p)=0$ otherwise.  In the notation of \S \ref{JP_elliott_03_sec_first_waystation}, with $\mathcal{J}$ the single real nonprincipal character $\modd{D}$ guaranteed by an application of Theorem \ref{JP_elliott_03_thm_B} to the function $g\mu$ with $r=2$, $2\sqrt{2}\varepsilon=1$, $c=1$, $\eta = 2^{-9}$, or the empty set if there is none such, the convolution factorisation $g\Lambda = g\mu * g\log$ affords the representation
\begin{align}
Y(g\Lambda,a,x) & = \sum\limits_{m \le y} g(m)\log m \, Y(g\mu, a\overline{m}, xm^{-1} ) \notag \\
& + \sum\limits_{r \le xy^{-1}} g(r)\mu(r) \left( Y(g\log, a\overline{r}, xr^{-1}) - Y(g\log, a\overline{r}, y) \right). \notag
\end{align}

We assume that $D^5 \le y \le x^{1/2}$.  The sum over $m$ is then
\begin{equation}
\ll \sum\limits_{m \le y} \frac{g(m) \log m}{m} \frac{x}{\varphi(D)\log x} \prod\limits_{D < p \le x} \left(1+ \frac{1}{p} \right)\left( \frac{\log D}{\log x} \right)^\eta \notag
\end{equation}
and, in view of the elementary bounds
\begin{equation}
\sum\limits_{m \le y} \frac{g(m)}{m} \ll \prod\limits_{D < p \le y} \left( 1+ \frac{1}{p} \right) \ll \frac{\log y}{\log D}, \notag
\end{equation}
at most a constant multiple of $\varphi(D)^{-1} x(\log y/\log D)^2 (\log D/\log x)^\eta$.

The treatment of the sum over $r$ is largely an exercise in the application of a sieve.  Since a version of it in some detail may be found in the appendix to \cite{elliottMFoAP6}, we confine ourselves to main points.

If $w \ge D^5$, $(a,D)=1$, then an application of the \emph{fundamental lemma} version of Selberg's sieve, c.f. \cite{Elliott1979}, Chapter 2, provides an estimate
\begin{equation}
\sum\limits_{\substack{m \le w \\ m \equiv a \modd{D}}} g(m) - \frac{1}{\varphi(D)} \sum\limits_{m \le w} g(m) \ll \frac{w}{\varphi(D)\log D} \exp\left(-\frac{2\log w}{\log D} \right), \notag
\end{equation}
the implied constant absolute.  An integration by parts yields
\begin{equation}
\sum\limits_{\substack{ y < m \le w \\ m \equiv a \modd{D}}} g(m)\log m - \frac{1}{\varphi(D)} \sum\limits_{y < m \le w} g(m)\log m \ll \frac{w}{\varphi(D)} \exp\left( - \frac{\log w}{\log D} \right), \notag
\end{equation}
uniformly for $D^2 \le y \le w$.

By partitioning the variable $m$ into residue classes $\modd{D}$, a similar upper bound holds for the sum
\begin{equation}
\varphi(D)^{-1} \sum\limits_{y < m \le w} g(m)\chi(m)\log m. \notag
\end{equation}
The initial sum over $r$ is thus
\begin{equation}
\ll \sum\limits_{r \le x/y} \frac{g(r)|\mu(r)|}{r}  \frac{x}{\varphi(D)} \exp\left( -\frac{\log(x/r)}{\log D} \right) \notag
\end{equation}
and, with $\theta = 1/\log D$,
\begin{equation}
\ll \frac{x}{\varphi(D)} \exp\left( - \frac{\log x}{\log D} \right) \sum\limits_{r \le x/y} \frac{g(r)}{r^{1-\theta}}. \notag
\end{equation}
An integration by parts combined with the estimate(s)
\begin{equation}
\sum\limits_{m \le t} g(m) \le \sum\limits_{\substack{m \le t \\ p \mid m \implies p > D^{1/2}}} 1 \ll t(\log D)^{-1}, \notag
\end{equation}
valid uniformly for $t \ge D$, shows the final sum over $r$ to be $\ll \exp\left( \log(x/y)(\log D)^{-1} \right)$.

Altogether
\begin{equation}
x^{-1}\varphi(D) Y(g\Lambda,a,y) \ll \left( \frac{\log y}{\log D} \right)^2 \left( \frac{\log D}{\log x} \right)^\eta + \exp\left(-\frac{\log y}{\log D} \right). \notag
\end{equation}
We choose $\log y = \log D(\log x/\log D)^{\eta/4}$.

The factor $g$ may be removed from $Y(g\Lambda,a,y)$ at an expense of $O(x^{1/2})$, negligible in comparison with the target error, and Lemma \ref{JP_elliott_03_lem_BB13} is established.

\noindent \emph{Proof of Linnik's theorem}.  We apply the corollary to Lemma \ref{JP_elliott_03_lem_BB13}.  According to Lemma \ref{JP_elliott_03_lem_I3}, the sum $\sum_{N^\gamma < p \le N} \chi(p)p^{-1}$ is bounded above by an absolute constant $c_1$, uniformly for $D \le N^\gamma$.  If $\chi(a) = -1$, then
\begin{equation}
\sum\limits_{N^\gamma < p \le N} p^{-1} (1-\chi(p)) > -\log \gamma - c_1 + O\left( \left( \frac{\log D}{\log N} \right)^\tau \right) \notag
\end{equation}
from which, with $\gamma \exp(2c_1) = 1$, Linnik's theorem follows at once.

The cases when $\chi(a) = 1$ are supplied by the following result, an elementary proof of which may be found as Lemma 13 in \cite{elliott2002millenium}.  For $0 < \beta < 1$, $y > 0$, and $\chi$ a real nonprincipal character $\modd{D}$, define
\begin{equation}
M_\beta(y) = \sum\limits_{\substack{ y^\beta < p \le y \\ \chi(p)=1}} p^{-1}. \notag
\end{equation}

\makeatletter{}

\begin{lemma} \label{JP_elliott_03_lem_BB14}
Let $96\beta < 1$, $\varepsilon > 0$, $y \ge 4$, $\chi(a)=1$.  Then
\begin{equation}
\sum\limits_{\substack{ y^\beta < p \le y \\ p \equiv a \modd{D}}} 1  = \left( 1 + O\left(\exp\left(-\frac{1}{7\beta}\log \frac{1}{\beta} \right) \right) + O(M_\beta(y)) \right) \frac{y}{D} A(D)+ O\left(y^{7/8+\varepsilon} \right), \notag
\end{equation}
where
\begin{equation}
A(D) = L(1,\chi)\prod\limits_{\substack{p \le y^\beta \\ \chi(p)=1}} \left(1-\frac{2}{p+2} \right)  \prod\limits_{\chi(p)=1} \left(1- \frac{3}{p^2} + \frac{2}{p^3} \right)\prod\limits_{\chi(p)=-1} \left(1-\frac{1}{p^2} \right), \notag
\end{equation}
the implied constant depending at most upon $\varepsilon$.
\end{lemma} 

With the roles of $\beta$, $y$ played by $\gamma$, $N$ respectively, for a sufficiently small absolute value of $\gamma$, either there is a prime $p \equiv a \modd{D}$ in the interval $(N^\gamma,N]$, or $M_\gamma(N)$ exceeds a further positive absolute constant $c_2$.  In the latter case
\begin{equation}
\sum\limits_{N^\gamma < p \le N} p^{-1} \left( 1+\chi(p) \right) \ge 2M_\gamma(N) \ge 2c_2 \notag
\end{equation}
and Linnik's theorem is again evident.
\end{example}

References in the following comments are confined to works narrowly connected to the present paper.

\emph{Comments}.  A systematic study of multiplicative functions with values in the complex unit disc, initiated by Delange in 1961 \cite{delange1961surlesfonctions}, received strong impulses from Wirsing \cite{wirsing1967}, 1967, and Hal\'asz \cite{halasz1968}, 1968.  Although the identity $\log n = \sum_{d \mid n} \Lambda(d)$ is employed in the works of Chebyshev, Wirsing seems to have been the first to apply it systematically to the study of general multiplicative functions.  The influence of these authors in the present paper is everywhere evident.

A version of Theorem \ref{JP_elliott_03_thm_01}, valid for $\alpha < 1/4$, may be found in Elliott \cite{elliottMFoAP7}.

An early version of Theorem \ref{JP_elliott_03_thm_01}, valid for $\alpha < 1$, under the slightly simplifying assumption that the modulus $D$ exceed an arbitrarily small power of $\log x$, is carried out in careful detail in the second author's Ph.D. thesis \cite{kish2013phd}.  Whilst the general outline follows that in Elliott \cite{elliottMFoAP6}, there are serious simplifying and improving modifications and a sharpened error term.  In particular, the step from an $L^2$ to an $L^\infty$ estimate is effected by the interpolation of a logarithm rather than the application of an approximate sieve identity.

Besides its generality, an important feature of Theorem \ref{JP_elliott_03_thm_01} is the uniformity in $D$, sufficient to establish Linnik's theorem on the size of the least prime in an arithmetic progression.  Indeed, within the same generality, this was a calibrating target in the series of papers \cite{elliottMFoAP1, elliottMFoAP2, elliottMFoAP3, elliottMFoAP4, elliottMFoAP5, elliottMFoAP6, elliottMFoAP7} by the first author, key arguments from which are subsumed in the present paper.  A proof of Linnik's theorem was achieved in the Illinois Millenial Conference paper \cite{elliott2002millenium}.

The recovery of a multiplicative function from the space generated by the one-parameter group $S_\tau$, $\tau \in \mathbb{R}$, is aided by the following result.

\makeatletter{}

\begin{lemma} \label{JP_elliott_03_lem_one_parameter_recovery}
The arithmetic function
\begin{equation}
\psi(n) = \sum\limits_{j=1}^k c_j n^{i\tau_j}, \quad n \ge 1, \ c_j \in \mathbb{C}, \ \tau_j \in \mathbb{R}, \notag
\end{equation}
is nontrivial and multiplicative if and only if $k=1$, $c_1 = 1$.
\end{lemma} 

\makeatletter{}

\noindent \emph{Proof of Lemma \ref{JP_elliott_03_lem_one_parameter_recovery}}.  For each prime $p$, $(m,p)=1$, multiplicativity ensures that
\begin{equation}
\psi(p)\psi(m) = \sum\limits_{j=1}^k c_j p^{i\tau_j} m^{i\tau_j}, \notag
\end{equation}
hence
\begin{equation}
\sum\limits_{j=1}^k c_j \left(p^{i\tau_j} - \psi(p) \right) m^{i\tau_j} = 0. \notag
\end{equation}
Then, if $\sigma > 1$,
\begin{align}
\sum\limits_{j=1}^k c_j & \left( p^{i\tau_j} -  \psi(p) \right) \sum\limits_{(m,p)=1} m^{-\sigma + i\tau_j}  \notag \\
& = \sigma \int_{1^-}^\infty y^{-\sigma-1} \sum\limits_{j=1}^k c_j \left( p^{i\tau_j} - \psi(p) \right) \sum\limits_{\substack{ m \le y \\ (m,p)=1}} m^{i\tau_j} \, dy = 0. \notag
\end{align}
Considering Euler products,
\begin{equation}
\sum\limits_{j=1}^k c_j \left( p^{i\tau_j} - \psi(p) \right)\left(1 - p^{-\sigma + i\tau_j} \right) \zeta(\sigma + i\tau_j) = 0, \quad \sigma > 1, \notag
\end{equation}
and by analytic continuation, since $\zeta(s)$ has only one singularity, a simple pole at $s=1$,
\begin{equation}
\sum\limits_{j=1}^k c_j \left(p^{i\tau_j} - \psi(p) \right) \left(1- p^{-z + i\tau_j} \right)\zeta(z+i\tau_j) = 0 \notag
\end{equation}
for all $z$ such that no $z+i\tau_j$, $j = 1, \dots, k$, has the value 1.  Allowing $z$ to approach each point $1-i\tau_j$ in turn, we see that
\begin{equation}
c_j \left(p^{i\tau_j} - \psi(p) \right) = 0, \quad 1 \le j \le k. \notag
\end{equation}

If $k \ge 2$, then $\psi(p) = p^{i\tau_1} = p^{i\tau_2}$, $\tau_1 \ne \tau_2$, say, i.e. there exists $\tau \ne 0$, real, such that $p^{i\tau} = 1$ for all primes $p$.  Then for $\sigma > 1$,
\begin{equation}
\zeta(\sigma) = \prod_p (1-p^{-\sigma})^{-1} = \prod_p (1-p^{-\sigma-i\tau})^{-1} = \zeta(\sigma + i\tau), \notag
\end{equation}
$\zeta(z) = \zeta(z+i\tau)$ provided $z \ne 1, 1-i\tau$.  This is untenable; $k=1$.

At this stage $\psi(n)$ can only be nontrivial and multiplicative if, for distinct primes $p$, $q$
\begin{equation}
c_1^2 p^{i\tau_1}q^{i\tau_1} = \psi(p)\psi(q) = \psi(pq) = c_1(pq)^{i\tau_1}, \notag
\end{equation}
$c_1^2 = c_1$, $c_1 = 1$.

This completes the proof.

An early realisation of Linnik's Large Sieve as an inequality attached to the action of a self-adjoint operator, and that such inequalities come in pairs, may be found in Elliott \cite{elliott1971inequalitieslargesieve}.

Functional analysis may also be applied to the study of complex-valued multiplicative functions with values outside the unit disc, c.f. Elliott \cite{Elliott1997}, \cite{elliott2010operatornorms}.  To offset the lack of an obvious referent for size, the function is compared to its absolute self.  We may continue to view nonvanishing completely multiplicative functions as characters on $\mathbb{Q}^*$ provided we abandon the requirement that the associated group representations be unitary.

A particular one-sided inequality related to Lemma \ref{JP_elliott_03_lem_I2} was privately circulated at the American Mathematical Society Research Community meeting on `The Pretentious View of Analytic Number Theory' held in Snowbird, Utah, Summer of 2011, attended by the second author.  The present two-sided, more general inequality, a slightly modified version of the theorem in the authors' paper \cite{elliottandkishone2012}, where a number of relevant further comments may be found, has a different proof.  Currently, all inequalities of related type rest upon a version of Lemma \ref{JP_elliott_03_lem_I3} developed from that of Elliott \cite{elliottMFoAP6}, with a variant argument in Elliott \cite{elliottMFoAP7}.

The results of \S \ref{JP_elliott_03_sec_taxonomy} elaborate the taxonomy of exceptional characters in terms of the support of the function $g$ carried out in the authors' paper \cite{elliottandkishtwo2012}.  A detailed study related to Lemmas \ref{JP_elliott_03_lem_A3} and \ref{JP_elliott_03_lem_A4} of the present work may be found in \S 3 of Elliott \cite{elliott1988additivearithfuncs}.

As was demonstrated in the Illinois Millenial paper \cite{elliott2002millenium}, to establish Linnik's theorem requires only a version of Theorem \ref{JP_elliott_03_thm_01} for functions supported on the primes in the interval $(D,x]$, and that at most one Dirichlet character $\modd{D}$ can be near to a given multiplicative function.

Lemma \ref{JP_elliott_03_lem_BB14}, adapted from a result of Heath-Brown \cite{heathbrown1990siegelzeros}, is derived by applying a sieve to an asymptotic estimate for the mean-value of the arithmetic function $n \mapsto \sum_{d \mid n} \chi(d)$ over an appropriate arithmetic progression.  It may be viewed as a localised descendent of the original argument employed by Dirichlet to establish the nonvanishing of $L(1,\chi)$ for a real character, $\chi$.

As a sampling of examples in the application of harmonic analysis on $\mathbb{Q}^*$ to problems in analytic number theory, and that apply methods related to those in the present paper:

Let $f_1$, $f_2$ be real-valued additive arithmetic functions.  Necessary and sufficient conditions for the arithmetic function $f_1(an+b) + f_2(An+B)$, with integers $a>0$, $A>0$, $b$, $B$, $aB \ne Ab$, to satisfy a weak limit law when suitably renormalised, obtained via the study of characters on $\left(\mathbb{Q}^* \times \mathbb{Q}^*\right)/\Gamma_k$, where $\Gamma_k$, subgroup of the direct product of two copies of $\mathbb{Q}^*$, is generated by elements of the form $(an+b) \times (An+B)$, $n \ge k \ge 1$, are established in the first author's memoir \cite{elliott1994correlationmemoir}; see also \cite{Elliott1985}.

A similar but more elaborate study of $f_1(n) + f_2(N-n)$, $1 \le n \le N$, is carried out in Elliott \cite{elliott2010valuedistribution}.

Structurally best possible bounds for the value concentration of an additive arithmetic function on integers of the form $p+1$, $p$ prime, or $N-p$, $p < N$, are established in Elliott \cite{elliott1994concentration}.

An upper bound on the size of the group $\mathbb{Q}^*/E_k$, where $E_k$ is the subgroup of $\mathbb{Q}^*$ generated by the shifted primes $p+1$, $p \ge k \ge 1$, is derived in Elliott \cite{elliott1995shiftedprimes1}.

A general discussion of factor groups of the direct product of finitely many copies of $\mathbb{Q}^*$, with attendant details of harmonic analysis may be found in Elliott \cite{elliott2002productrepresentations}; an exotic example in the application of such harmonic analysis to the product/quotient representations of positive rationals by products of variously shifted primes may be found in Elliott \cite{elliott2007kanazawa}.

\bibliographystyle{amsplain}
\bibliography{MathBib}

\providecommand{\bysame}{\leavevmode\hbox to3em{\hrulefill}\thinspace}
\providecommand{\MR}{\relax\ifhmode\unskip\space\fi MR }
\providecommand{\MRhref}[2]{%
  \href{http://www.ams.org/mathscinet-getitem?mr=#1}{#2}
}
\providecommand{\href}[2]{#2}
\begin{thebibliography}{10}

\bibitem{delange1961surlesfonctions}
H.~Delange, \emph{Sur les fonctions arithm\'etiques multiplicatives}, Ann. Sci.
  \'Ecole Norm. Sup. (3) \textbf{78} (1961), 273--304.

\bibitem{elliott1971inequalitieslargesieve}
P.~D. T.~A. Elliott, \emph{On inequalities of large sieve type}, Acta Arith.
  \textbf{18} (1971), no.~1, 405--422, Davenport Memorial Volume.

\bibitem{Elliott1979}
\bysame, \emph{Probabilistic number theory {I}: Mean-value theorems},
  Grundlehren der Mathematischen Wissenschaften, vol. 239, Springer-Verlag, New
  York, 1979.

\bibitem{Elliott1985}
\bysame, \emph{Arithmetic functions and integer products}, Grundlehren der
  Mathematischen Wissenschaften, vol. 272, Springer-Verlag, New York, 1985.

\bibitem{elliottMFoAP1}
\bysame, \emph{Multiplicative functions on arithmetic progressions},
  Mathematika \textbf{34} (1987), no.~2, 199--206.

\bibitem{elliott1988additivearithfuncs}
\bysame, \emph{Additive arithmetic functions on intervals}, Math. Proc. Camb.
  Phil. Soc. \textbf{103} (1988), no.~1, 163--179.

\bibitem{elliottMFoAP2}
\bysame, \emph{Multiplicative functions on arithmetic progressions. {II}},
  Mathematika \textbf{35} (1988), no.~1, 38--50.

\bibitem{elliottMFoAP3}
\bysame, \emph{Multiplicative functions on arithmetic progressions. {III}.
  {T}he large moduli}, A tribute to {P}aul {E}rd{\H o}s, Cambridge University
  Press, Cambridge, 1990, pp.~177--194.

\bibitem{elliottMFoAP4}
\bysame, \emph{Multiplicative functions on arithmetic progressions. {IV}. {T}he
  middle moduli}, J. London Math. Soc. (2) \textbf{41} (1990), no.~2, 201--216.

\bibitem{elliottMFoAP5}
\bysame, \emph{Multiplicative functions on arithmetic progressions. {V}.
  {C}omposite moduli}, J. London Math. Soc. (2) \textbf{41} (1990), no.~3,
  408--424.

\bibitem{elliottMFoAP6}
\bysame, \emph{Multiplicative functions on arithmetic progressions. {VI}.
  {M}ore middle moduli}, J. Number Theory \textbf{44} (1993), no.~2, 178--208.

\bibitem{elliott1994concentration}
\bysame, \emph{The concentration function of additive functions on shifted
  primes}, Acta Math. \textbf{173} (1994), no.~1, 1--35.

\bibitem{elliott1994correlationmemoir}
\bysame, \emph{On the correlation of multiplicative and the sum of additive
  arithmetic functions}, Mem. Amer. Math. Soc. \textbf{112} (1994), no.~538,
  viii+88.

\bibitem{elliott1995shiftedprimes1}
\bysame, \emph{The multiplicative group of rationals generated by the shifted
  primes, {I}}, J. Reine Angew. Math. \textbf{463} (1995), 169--216.

\bibitem{Elliott1997}
\bysame, \emph{Duality in analytic number theory}, Cambridge Tracts in
  Mathematics, vol. 122, Cambridge University Press, Cambridge, 1997.

\bibitem{elliott2002millenium}
\bysame, \emph{The least prime primitive root and {L}innik's theorem}, Number
  theory for the millennium {I}, Proceedings of the Millennial Conference on
  Number Theory, University of Illinois at Urbana-Champaign, May 21--26, 2000
  (M.A. Bennett, B.C. Berndt, N.~Boston, H.G. Diamond, A.J. Hildebrand, and
  W.~Philipp, eds.), A K Peters, Natick, 2002, pp.~393--418.

\bibitem{elliottMFoAP7}
\bysame, \emph{Multiplicative functions on arithmetic progressions. {VII}.
  {L}arge moduli}, J. London Math. Soc. (2) \textbf{66} (2002), no.~1, 14--28.

\bibitem{elliott2002productrepresentations}
\bysame, \emph{Product representations by rationals}, Number Theoretic Methods:
  Future Trends, Proceedings of the Second China-Japan Seminar, Iizuka, Japan,
  March 12--16, 2001 (Shigeru Kanemitsu and Chaohua Jia, eds.), Dev. Math.,
  vol.~8, Kluwer Acad. Publ., Dordrecht, 2002, pp.~119--150.

\bibitem{elliott2007kanazawa}
\bysame, \emph{The ramifications of a shift by 2}, Probability and number
  theory, Proceedings of the International Conference on Probability and Number
  Theory, Kanazawa, Japan, June 20--24, 2005 (Shigeki Akiyama, Kohji Matsumoto,
  Leo Murata, and Hiroshi Sugita, eds.), Adv. Stud. Pure Math., vol.~49, Math.
  Soc. Japan, Tokyo, 2007, pp.~69--77.

\bibitem{elliott2010valuedistribution}
\bysame, \emph{The value distribution of additive arithmetic functions on a
  line}, J. Reine Angew. Math. \textbf{642} (2010), 57--108.

\bibitem{elliott2010operatornorms}
\bysame, \emph{Operator norms and the mean-values of multiplicative functions},
  Functions in Number Theory and their Probabilistic Aspects, Proceedings of
  the International Conference on Functions in Number Theory and their
  Probabilistic Aspects, Research Institute for Mathematical Sciences, Kyoto
  University, Japan, December 13--17, 2010 (Kohji Matsumoto~(Editor in~Chief),
  Shigeki Akiyama, Katusi Fukuyama, Hitoshi Nakada, Hiroshi Sugita, and Akio
  Tamagawa, eds.), vol. B34, RIMS Kokyuroku Bessatsu, 2012, pp.~81--102.

\bibitem{elliottandkishone2012}
P.~D. T.~A. Elliott and J.~Kish, \emph{A large sieve inequality for {E}uler
  products}, arXiv:1203.0804 (2012).

\bibitem{elliottandkishtwo2012}
\bysame, \emph{Multiplicative functions and a taxonomy of {D}irichlet
  characters}, arXiv:1208.0051 (2012).

\bibitem{halasz1968}
G.~Hal{\'a}sz, \emph{\"{U}ber die {M}ittelwerte multiplikativer
  zahlentheoretischer {F}unktionen}, Acta Math. Acad. Sci. Hungar. \textbf{19}
  (1968), 365--403.

\bibitem{heathbrown1990siegelzeros}
D.~R. Heath-Brown, \emph{Siegel zeros and the least prime in an arithmetic
  progression}, Quart. J. Math. Oxford Ser. (2) \textbf{41} (1990), no.~164,
  405--418.

\bibitem{kish2013phd}
J.~Kish, \emph{Harmonic analysis on the positive rationals: Multiplicative
  functions and exceptional {D}irichlet characters}, Ph.D. thesis, University
  of Colorado Boulder, 2013.

\bibitem{shiu1980bruntitchmarsh}
P.~Shiu, \emph{A {B}run-{T}itchmarsh theorem for multiplicative functions}, J.
  Reine Angew. Math. \textbf{313} (1980), 161--170.

\bibitem{titchmarsh1986zeta}
E.C. Titchmarsh, \emph{The theory of the {R}iemann zeta-function}, second ed.,
  Oxford University Press, New York, 1986, Revised by D.R. Heath-Brown.

\bibitem{wirsing1967}
E.~Wirsing, \emph{Das asymptotische {V}erhalten von {S}ummen \"uber
  multiplikative {F}unktionen. {II}}, Acta Math. Acad. Sci. Hungar. \textbf{18}
  (1967), 411--467.

\end{thebibliography}

First published by transmission September 4, 2013.  An occasional typographical oversight may remain.

\end{document}